\documentclass[11pt]{article}

\usepackage[top=2cm, bottom=3cm, left=2cm, right=2cm]{geometry}

\usepackage[abbrev]{amsrefs}

\makeatletter
\long\def\unmarkedfootnote#1{{\long\def\@makefntext##1{##1}\footnotetext{#1}}}
\makeatother

\usepackage{graphicx}
\usepackage{amsmath, amssymb,
fancybox,enumerate,color,amsthm}

\begin{document}

\renewcommand{\theequation}{\arabic{section}.\arabic{equation}}

\newtheorem{theorem}{\bf Theorem}[section]
\newtheorem{lemma}{\bf Lemma}[section]
\newtheorem{proposition}{\bf Proposition}[section]
\newtheorem{corollary}{\bf Corollary}[section]
\newtheorem{remark}{\bf Remark}[section]
\newtheorem{example}{\bf Example}[section]
\newtheorem{definition}{\bf Definition}[section]

\newcommand{\supp}{\mathop{\mathrm{supp}}}
\newcommand{\PA}{\partial}
\newcommand{\ve}{\varepsilon}
\newcommand{\TR}{\textcolor{red}}
\newcommand{\TB}{\textcolor{blue}}
\newcommand{\ds}\displaystyle

\newcommand{\re}{{\Bbb R}}
\newcommand{\nid}{\noindent}

\newcommand{\nn}{\nonumber}
\newcommand{\hu}{\hat{u}}
\newcommand{\tw}{\widetilde{w}}
\newcommand{\D}{D_0^{1,2}}
\newcommand{\DK}{D_{0,K}^{1,2}}
\newcommand{\rhp}{\rightharpoonup}
\def\R{{\mathbb R}}
\newcommand{\la}{\langle}
\newcommand{\ra}{\rangle}
\newcommand{\ls}{\{}
\newcommand{\rs}{\}_{n=1,2\cdots}}
\newcommand{\hv}{\hat{v}}
\newcommand{\hvp}{\widehat{\vp}}
\newcommand{\hw}{\widehat{w}}
\newcommand{\cv}{\bar{v}}
\newenvironment{prf}[1]
   {{\noindent \bf Proof of {#1}.}}{\hfill \qed}
\newcommand{\eq}[1]{{\begin{equation}#1\end{equation}}}

\makeatletter
 \@addtoreset{equation}{section}
\makeatother
\def\theequation{\arabic{section}.\arabic{equation}}

\title{
Singular solutions of semilinear elliptic equations with \\ exponential nonlinearities in 2-dimensions
}
\author{\small
Yohei Fujishima\\
\small
Department of Mathematical and Systems Engineering,
Faculty of Engineering,
\\
\small
Shizuoka University, \\
\small
3-5-1 Johoku, Hamamatsu 432-8561, Japan, 
\\
\small
{fujishima@shizuoka.ac.jp}
\vspace{5pt}\\
\quad\\
\small
Norisuke Ioku
\\
\small
Mathematical Institute, 
\\
\small
Tohoku University, \\
\small
Aramaki 6-3, Sendai 980-8578, Japan,
\\
\small
{ioku@tohoku.ac.jp}
\vspace{5pt}\\
\quad\\
\small
Bernhard Ruf
\\
\small
Accademia di Scienze e Lettere - Istituto Lombardo,\\
\small
via Borgonuovo, 25, Milano 20121, Italy,
\\
\small
{bruf001@gmail.com}
\\
\quad\\
{\small and}\\
\quad\\
\small
Elide Terraneo\\
\small
Dipartimento di Matematica ``F. Enriques'',\\
\small
Universit\`a degli Studi di Milano,\\
\small
via C. Saldini 50, Milano 20133, Italy,
\\
\small
{elide.terraneo@unimi.it}
}

\date{
}
\maketitle

\abstract{
By introducing a new classification of the growth rate of exponential functions, singular solutions for 
$-\Delta u=f(u)$ in 2-dimensions with exponential nonlinearities are constructed.
The strategy is to introduce a model nonlinearity ``close" to $f$, which admits an explicit singular solution. Then, using a transformation as in \cite{FI}, one obtains an approximate singular solution, and then one concludes by a suitable fixed point argument. 
Our method covers a wide class of nonlinearities in a unified way, e.g., 
$f(u) = u^re^{u^q}\ (q>1,r\in \mathbb{R}), f(u) = e^{u^{q}+u^r}\ (q>1,\ q/2>r>0)$ or $f(u) = e^{e^u}$.
As a special case,
our result contains a pioneering contribution by Ibrahim--Kikuchi--Nakanishi--Wei~\cite{IKNW} for $u(e^{u^2}-1)$.
}

\ \\
\noindent
{\small
{\bf Keywords}: 
Semilinear elliptic equations, singular solutions, 2-dimensions, exponential nonlinearities
\vspace{5pt}
\newline
{\bf 2010 MSC}: Primary: 35J61; Secondly: 35B40, 35A24 
\vspace{5pt}
}

\newpage
\section{Introduction}

We are concerned with singular radial solutions of the semilinear elliptic equation
\begin{equation}\label{eq:1.1}
-\Delta u = f(u)
\end{equation}
in a punctured ball $B_R(0)\setminus \{0\}\subset \mathbb{R}^N$,
where $N\ge 2$.
We say $u\geq 0$ is a singular solution of 
\eqref{eq:1.1}
if there exists $r_0$ such that $u\in C^2(B_{r_0}\setminus \{0\})$ satisfies 
$-\Delta u = f(u)$ in $B_{r_0}\setminus \{0\}$ in the classical sense
and $u(x)\to \infty$ as $|x|\to 0$. For an overview of such results, we refer to the books by Quittner--Souplet \cite{QS} and by V\'eron \cite{V}.

Let us recall some existence and uniqueness results of singular solutions 
of the equation \eqref{eq:1.1} for the particular case $f(u)=u^p$ for $N\ge 3$, i.e.,
\begin{equation}\label{eq:poli}
-\Delta u =u^p.
\end{equation}
When $1<p<\frac{N}{N-2}$, 
Lions~\cite{L}
showed that 
the equation~\eqref{eq:poli} has a singular solution which blows up with the same speed as the fundamental solution $|\cdot|^{2-N}$
of the Laplace equation.
For $\frac{N}{N-2}\leq p<\frac{N+2}{N-2}$, Ni--Sacks~\cite{NS} (see also Aviles~\cite{A} and Chen--Lin~\cite{CL}) proved the existence of infinitely many singular solutions of \eqref{eq:poli}. Moreover,   $u(x) \sim |x|^{2-N}(-\log|x|)^{\frac{2-N}2}$ near $0$, if $p=\frac{N}{N-2}$, and  $u(x) \sim |x|^{-\frac{2}{p-1}}$ near $0$, if $\frac{N}{N-2}<p<\frac{N+2}{N-2}$. 
These singular solutions also satisfy the equation~\eqref{eq:poli} in a non-punctured ball in the sense of distributions.
Finally, in the case $p >\frac{N}{N-2}$, it is easy to see that
\eqref{eq:poli} has the exact singular solution 
\begin{equation}\label{eq:explicit}
v(x)=L_p\, |x|^{-\frac{2}{p-1}},\qquad 
L_p
:=
\Big\{
\frac{2}{p-1}\big(N-2-\frac{2}{p-1}\big)
\Big\}^{\frac{1}{p-1}},
\end{equation}
which is  also a  distributional solution.
In addition, Caffarelli--Gidas--Spruck~\cite{CGS} showed that \eqref{eq:poli} has for $p=\frac{N+2}{N-2}$ a continuum of singular solutions. 
On the other hand, if $p> \frac{N+2}{N-2}$, it was proved by Serrin--Zou~\cite{SZ} that the radial singular solution of \eqref{eq:poli}  is unique.
As a corresponding case to $p=\infty$
for $N\ge 3$,
exponential type nonlinearities 
were investigated and singular solutions to the equation \eqref{eq:1.1} were constructed by Mignot--Puel \cite{MP} for $f(u)=e^u$,
Kikuchi--Wei~\cite{KW} for $f(u)=e^{u^q} (q>0)$,
Ghergu--Goubet~\cite{GG} for multiple exponential nonlinearities.
More general nonlinearities
in the form of $f(u)=u^p+\sigma(u)$
or $f(u)=e^u+\sigma(u)$ with a lower order term $\sigma(u)$
have been widely studied. 
Since there are a vast number of researches, we only refer to
\cite{J,M2,M3,MN2,MN3}; see also references therein.

\par \bigskip

In contrast with the case $N\ge 3$, much less is known for the two dimensional case $N=2$.
Lions~\cite{L} proved that a similar result to $N\ge 3$ holds for 
the
power type nonlinearity $f(u)=u^p$, namely,
the equation~\eqref{eq:poli} in 2 dimensions has a singular solution which
blows up with the same speed as the fundamental solution $-\log |x|$ of the Laplace equation.
For $f(u)=e^u$, Tello~\cite{T} showed that all solutions are given by an explicit two-parameter family of solutions, containing singular solutions  with growth like $-\log |x|$ near the origin, and containing also a one parameter family of regular solutions (see also Taliaferro\cite{Ta1,Ta2}). 
For functions $f$ with higher growth than exponential, we mention recent results for $f(s) \sim e^{s^2}$. In \cite{FR} it was noted that 
\begin{equation}\label{exp} 
	u(r) = (-2\log r)^{1/2} \ \hbox{ solves the equation  \eqref{eq:1.1} for }\ f(s) = \frac 1{s^3}\ e^{s^2}.
\end{equation}
Building on this solution, it was shown in \cite{IRT1} that there exists $R > 0$ such that the boundary value problem
\begin{equation}\label{Mpl}\left\{
\begin{array}{ll}  -\Delta u = g(u) \  \ \hbox{in} \ B_R\setminus\{0\},
	\vspace{0.2cm}\\
	\ \  \ \    \ u = 0 \ \ \ \quad  \hbox{on} \ \partial B_R,
\end{array} \right.
\end{equation}
has a singular radial and positive solution with growth $(-2\log|x|)^{1/2}$ near $0$, for the nonlinearity
$$
g(s) = \left\{ \begin{array}{ll} e^{s^2}{/s^3},\  \ s\ge b, \vspace{0.2cm}\\
                as^2,	\  \ 0 \le s \le b,
 \end{array}\right. \quad  \text{where}\ a, b\ \hbox{are constants such that }\ g \in C^1\bigl(
[0,\infty) \bigr).
 $$
This solution is also a distributional solution.
Independently,  Ibrahim--Kikuchi--Nakanishi--Wei~\cite{IKNW} constructed  a singular distributional solution for 
\eqref{Mpl}
with the  nonlinearity $g(u)=u(e^{u^2}-1)$.
The asymptotic behavior of this solution at the origin is $\big(-2\log|x| - 2\log(-\log|x|)-2\log2\big)^{1/2}$.
\par \medskip
Inspired by the solution \eqref{exp}, 
one observes (see also \cite
{GGP}) that there exist explicit 
singular solutions of $-\Delta u = g(u)$  
for certain specific nonlinearities
 with growth $g(s) \sim e^{s^q}, q > 1$. 
Indeed, for 
$$ g(s) = \frac 4{q\,q'}\frac{e^{s^q}}{s^{2q-1}}\ , \ q > 1 \ , \ \hbox{ the solutions are} \quad  v(x) = (-2\log|x|)^{1/q},
$$
while for 
$$ g(s) = 4\, \frac {e^{e^s}}{e^{2s}} \ , \ \hbox{ the solution is} \quad v(x) = \log(-2\log|x|).
$$
These are very specific nonlinearities. The aim of the present article is to give generalizations to a wide class of nonlinearities of this type. 

\par \bigskip
Recently,
a possible direction to treat general nonlinearities for $N \ge 3$ was proposed 
in \cite{FI}. Indeed, in \cite{FI} it was proved  that for any couple of positive and increasing functions $f$, $g$ such that
\begin{equation*}
f, \ g \in C^1([0,\infty)), \ \ \   {\rm and}\ \ \
 F(s)=\int_s^{\infty}\frac1{f(\tau)}\ d\tau<\infty, \  \ \  G(s)=\int_s^{\infty}\frac1{g(\tau)}\ d\tau<\infty
\end{equation*}
if the function $v$
 satisfies 
\[
-\Delta v=g(v),
\] 
then
\begin{equation}\label{eq:u}
	\tilde u:=F^{-1}\Bigl[G[v(x)]\Bigr] \ 
\end{equation}
satisfies
\begin{equation}\label{eq:transformation}
	-\Delta \tilde u=f(\tilde u)
	+\frac{|\nabla{\tilde u}|^2}{f(\tilde u)F(\tilde u)}\Big[ g'(v)G(v)-f'(\tilde u)F(\tilde u)\Big].
\end{equation}
Therefore, thanks to 
the transformation~\eqref{eq:u} and
the transformed equation~\eqref{eq:transformation},
one can connect a general nonlinearity $f$
to a {\it model nonlinearity} $g$. More precisely  the strategy is the following. First, one proceeds by classifying    general nonlinearities $f$ in terms of their {\it growth rate}, measured by the  H\"older conjugate  of the exponent  
\begin{equation}\label{eq:A}
	A := \lim_{s\to \infty}f'(s)F(s)\ .
\end{equation}
Secondly, one can  connect the nonlinearity $f$ to the model nonlinearity $g(s)=s^{A'}$, where $A'$ is the  H\"older conjugate  of the exponent  $A$, if $A>1$ and to $g(s)={e}^s$, if $A=1$.
We remark that for the model nonlinear terms $g(s)$, it is easy to see that  $g'(s)G(s)\equiv A$. 
Then, if $A>1$ 
the function $v(x)$ defined in \eqref{eq:explicit} with $p=A'$
satisfies $-\Delta v=v^{A'}$, and we have
\[
G\bigl(v(x)\bigr)=\int_v^{\infty}\frac{1}{s^{A'}}ds=(A-1)v^{-\frac{1}{A-1}}=(A-1)\Bigl(L_{A'}|x|^{-\frac{2}{A'-1}}\Bigr)^{-\frac{1}{A-1}}
=\frac{|x|^2}{2N-4A}.
\]
In conclusion, 
one can expect that 
\[
\tilde u:=F^{-1}\Bigl[G[v(x)]\Bigr]
=F^{-1}\Bigl[\frac{|x|^2}{2N-4A}\Bigr]
\]
is an approximate singular solution of $-\Delta u=f(u)$
for general nonlinearities which satisfy
$A=\lim_{s\to \infty}f'(s)F(s)$.
Indeed, 
Miyamoto~\cite{M} constructed a singular solution of \eqref{eq:1.1} with $f\in C^2([0,\infty))$ satisfying
$A<\frac{N+2}{4}$ in the form 
\[
u(x)=F^{-1}\Big[\frac{|x|^2}{2N-4A}(1+o(1)\Big]\ \ \text{as}\ \ |x|\to 0.
\]
Later, Miyamoto--Naito~\cite{MN} proved 
that the singular solution constructed in \cite{M} is 
the unique radially symmetric singular solution of \eqref{eq:1.1}.
Since the condition $A<\frac{N+2}{4}$ is equivalent to 
$p>\frac{N+2}{N-2}$ for $f(u)=u^p$, the uniqueness result \cite{MN} is a natural generalization of \cite{SZ}.
Their uniqueness result 
is applicable for wide classes of nonlinearities.
For example, $f(u)=u^p$ with $p>\frac{N+2}{N-2}$, 
$f(u)=u^p+u^q\ (p>\frac{N+2}{N-2},p>q>0)$
$f(u)=u^re^{u^q}\ (r\in \mathbb{R},q>0)$, and multiple exponential functions can be treated. 
Further application of the exponent $A$ and the transformation~{\eqref{eq:u}}
can be found in \cite{FHIL,FI2,FI3,FI4,IRT1}.

\bigskip
\section{Main results}

In this paper, we construct a singular solution of the equation \eqref{eq:1.1} with $N = 2$ for general exponential type nonlinearities in terms of the classification of nonlinearities by the exponent $A$ and the transformation~{\eqref{eq:u}}.
Since $A=1$ for all the above mentioned exponential nonlinearities,
a more precise classification adapted to exponential type nonlinearities is  required.

\par \medskip
To this end, we define 
\[
\frac{1}{B_1[f](s)}
:= 
\bigl(-\log F(s)\bigr)\big[1-f'(s)F(s)\big]
\]
and its de l'Hospital form by
\[
\frac{1}{B_2[f](s)}
:=
\frac{\Bigl(1-f'(s)F(s)\Bigr)'}
{\left(\frac{1}{-\log F(s)}\right)'}
= 
f'(s)F(s)
\Bigl(-\log F(s)\Bigr)^2
\left[
\frac{f(s)f''(s)F(s)}{f'(s)}-1
\right].
\]
Then we define 
\begin{equation*}
\displaystyle B:=\lim_{s\to \infty}B_2[f](s).
\end{equation*}
It is easy to check that 
$g(s)=\frac{4}{qq'}\frac{e^{s^q}}{s^{2q-1}}\ (q>1)$
satisfies $B=q'$
and
$g(s)=4\frac{e^{e^s}}{e^{2s}}$
satisfies $B=1$.
More generally,
\[
f(s)=s^re^{s^q}\ (q>1,\ r\in \mathbb{R}),\quad e^{s^q+s^r}\ (q>1,\ q>r>0),\quad 
e^{s^q\big(\log (e+s)\bigr)^r}\ (q>1,r\in\mathbb{R})
\]
have $B=q'
$, 
and multiple exponential functions as
\[
f(s)
=
e^{e^{s^q}}\ (q\ge 1),\quad e^{e^{e^s}},\quad e^{e^{\cdots^{e^s}}}
\]
have $B=1$. 

Using the exponent $B$, we construct 
a singular solution to \eqref{eq:1.1} under the following general hypotheses:
\begin{itemize}
\item[($f1$)]
\label{f1}
$f\in 
 C^2\bigl((s_0,\infty)\bigr)$,
$f(s)>0$, $F(s)<\infty$, and $f'(s)>0$ for $s>s_0$
with some $s_0>0$,
where
\[
F(s):=\int_s^{\infty}\frac{1}{f(\tau)}d\tau.
\]
\item[($f2$)]
\label{f2}
The limits 
$\displaystyle \lim_{s\to \infty}f'(s)F(s)$
and 
$\displaystyle \lim_{s\to \infty}B_2[f](s)$
exist and satisfy
\begin{equation}
\label{eq:1.2cc}
\displaystyle \lim_{s\to \infty}f'(s)F(s)=1
\end{equation}
and
\begin{equation}
\label{eq:1.2c}
\displaystyle B:=\lim_{s\to \infty}B_2[f](s)\in [1,\infty).
\end{equation}
\end{itemize}
Fix $f$ satisfying $(f1)$ and $(f2)$.
Define $v$, $g$,  and $G$  by
\begin{equation}\label{eq:1.3a}
v(x):=
\Big(\log \frac{1}{|x|^2}\Big)^{\frac{1}{B'}},
\quad
g(s):=\frac{4}{BB'}s^{1-2B'}e^{s^{B'}},
\quad
G(s):=\int_s^{\infty}\frac{d\tau}{g(\tau)}=\frac{B}{4}\frac{s^{B'}+1}{e^{s^{B'}}}
\quad \text{if\ $B>1$},
\end{equation}
and
\begin{equation}\label{eq:1.3b}
v(x):=
\log\log\frac{1}{|x|^2},
\quad
g(s):=4\frac{e^{e^s}}{e^{2s}},
\quad
G(s):=\int_s^{\infty}\frac{d\tau}{g(\tau)}
=
\frac{1}{4}\frac{e^s+1}{e^{e^s}}
\quad \text{if\ $B=1$},
\end{equation}
as model nonlinearities and corresponding explicit solutions. 
Remark that
\begin{equation}\label{eq:1.4a}
G\bigl(v(x)\bigr)
=
\frac{B}{4}|x|^2\ \Big(\log\frac{1}{|x|^2}+1\Big)
\end{equation}
for $B\ge 1$.
Let 
\begin{equation}\label{eq:1.4}
\tilde u(x):=F^{-1}\big[G(v(x))\bigr].
\end{equation}
Our aim is to show that $\tilde u$ is an approximate solution of \eqref{eq:1.1}, and by estimating the second term of
the right hand side in
\eqref{eq:transformation} as an error, we will prove the existence of a singular solution.
To state our result, set
\begin{equation}\label{eq:1.5}
\begin{aligned}
R_1(x)
:=
\bigg|\frac{1}{B_1[f]\bigl(\tilde u(x)\bigr)}-\frac{1}{B_1[g]\bigl(v(x)\bigr)}\bigg|\ ,\quad 
R_2(x)
:=
\bigg|\frac{1}{B_2[f]\bigl(\tilde u(x)\bigr)}-\frac{1}{B_2[g]\bigl(v(x)\bigr)}\bigg|\ ,
\end{aligned}
\end{equation}
and make the  further assumption
\begin{equation}
\label{eq:1.1a}
\lim_{x \to 0}(-\log |x|)^{\frac{1}{2}}\left[R_1(x)+R_2(x)\right]=0.
\end{equation}

Our first statement shows the existence of a singular solution of $-\Delta u =f(u)$.

\begin{theorem}\label{theorem:1.2}
Assume that $f$ satisfies 
(f1), (f2), and 
\eqref{eq:1.1a},
where $g$ is the model nonlinearity as given in \eqref{eq:1.3a} or  \eqref{eq:1.3b} with corresponding solution $v$.
Then there exists a singular solution $u$ of $-\Delta u=f(u)$ (locally around the origin)
satisfying
\begin{equation}\label{eq:1.2e}
u(x)
=
\tilde u(x)
+
O\Big(f\bigl(\tilde u(x)\bigr)F\bigl(\tilde u(x)\bigr)\sup_{|y|\le |x|}
\bigl(R_1(y)+R_2(y)\bigr)\Big)
\qquad \text{as}\  |x|\to 0.
\end{equation}
\end{theorem}
\par \bigskip

We do not 
know if the hypothesis~\eqref{eq:1.1a} is necessary, but currently we are not able to avoid it. 
The following is a typical example to which  our result applies. 
See Section~\ref{section:4}
for further examples. 
\begin{corollary}\label{corollary:1.1}
For $q>1$ and $r\in \mathbb{R}$,
there exists a singular solution $u$ of $-\Delta u=u^re^{u^q}$
satisfying
\begin{equation}\label{eq:1.2}
u(x)
=
\bigg(
\log\frac{1}{|x|^2}-\frac{2q+r-1}{q}\log\log\frac{1}{|x|^2}+\log\frac{4(q-1)}{q^2}
\bigg)^{\frac{1}{q}}
+O\bigg(\frac{\log(-\log|x|)}{(-\log|x|)^{2-\frac{1}{q}}}\bigg)
\ \ \text{as $|x|\to 0$.}
\end{equation}
\end{corollary}

\begin{remark}
Ibrahim--Kikuchi--Nakanishi--Wei~\cite{IKNW}
constructed a singular solution for $f(u)=u(e^{u^2}-1)$.
Both the singularity and the remainder term of our solution \eqref{eq:1.2} coincide with theirs by taking $r=1$ and $q=2$ in Corollary~\ref{corollary:1.1}.
Moreover, 
if $q=2,r=-3$,
the singular solution constructed in Theorem~\ref{theorem:1.2} coincides with the explicit solution 
$v=\left(-2\log |x|\right)^{1/2}$ which is obtained 
 in~\cite{FR},
since it follows from $f(s)=s^{-3}e^{s^2}=g(s)$
that $\tilde u=v$ and $R_1=R_2=0$.
\end{remark}

\par \medskip

One can extend the singular solution constructed in Theorem~\ref{theorem:1.2} to a solution of the Dirichlet problem
by the shooting method.

\begin{theorem}\label{theorem:1.1}
Let $f$ satisfy
(f1), (f2), and \eqref{eq:1.1a} with respect to 
the model nonlinearity $g$.
Assume further that 
$f\in 
C^1\bigl([0,\infty)\bigr)$
and
$f(s)>0$, $F(s)<\infty$, $f'(s)>0$ for all $s>0$.
Then there exist $R>0$ and 
$u_{\infty}\in 
C^{2}(B_R\setminus\{0\})
$ 
satisfying
$\eqref{eq:1.2e}$
and
\begin{equation}\label{eq:1.1aaa}
\left\{
\begin{aligned}
-\Delta u
&
=f(u)
&& \text{in}\ B_R\setminus\{0\},
\\
u
&
=0
&& \text{on}\ \partial B_R,
\end{aligned}
\right.
\end{equation}
in the classical sense.
Furthermore,
$u_{\infty}$ satisfies
$-\Delta u_{\infty}=f(u_{\infty})$ in $B_R$
in the sense of distributions.
\end{theorem}
\par \medskip
Let us give further remarks about the exponents 
$A$ and $B$.
\par \medskip
\begin{remark}
It is pointed out in \cite{M} (see also \cite[Remark~1.2]{FI})
that
the exponent $A$ defined in \eqref{eq:A} coincides with
the limit 
$\displaystyle \lim_{s\to \infty}\frac{f'(s)^2}{f(s)f''(s)}$
which is introduced 
by Dupaigne--Farina~{\rm \cite{DF}}
as the generalized H\"older conjugate of the growth rate of $f$.
Indeed, by de l'Hospital rule, it holds
\[
 \frac{1}{A}=\lim_{s\to \infty}\dfrac{1/f'(s)}{F(s)}
 =\lim_{s\to \infty}\dfrac{\big(\frac{1}{f'(s)}\big)'}{\left(F(s)\right)'}
 =\lim_{s\to \infty}\frac{f(s)f''(s)}{f'(s)^2}
 =\lim_{s\to \infty}\bigg(1-\frac{\big(\frac{f(s)}{f'(s)}\big)'}{s'}\bigg)
 =\lim_{s\to \infty}\Big(1-\frac{f(s)}{sf'(s)}\Big),
\]
where the limit $\displaystyle \lim_{s\to \infty}\frac{sf'(s)}{f(s)}$ is regarded as the growth rate of a general nonlinear term $f$.
\end{remark}
\par \medskip
\begin{remark}\label{remark:1.3}
If (f1) and (f2) are satisfied then it holds
\begin{equation}\label{eq:1.2aaaa}
\lim_{s\to \infty}B_1[f](s)=B,
\end{equation}
since by 
 de l'Hospital rule
\begin{equation}\label{eq:1.11a}
\lim_{s\to \infty}\frac{1}{B_1[f](s)}
=
\lim_{s\to \infty}\frac{\Bigl(1-f'(s)F(s)\Bigr)'}{\Big(\frac{1}{-\log F(s)}\Big)'}
=
\lim_{s\to \infty}\frac{1}{B_2[f](s)}=\frac{1}{B}.
\end{equation}
Furthermore, it holds
\begin{equation}\label{eq:1.2b}
f'(s)F(s)\to 1^-\ \ \ {\rm  as}\ \  s\to \infty.
\end{equation}
Indeed, by \eqref{eq:1.2aaaa}
for any $\varepsilon\in \left(0,\frac{1}{B}\right)$
there exists $s_0$ such that
\[
\frac{1}{B}-\varepsilon 
<
\bigl(-\log F(s)\bigr)[1-f'(s)F(s)]
<
\frac{1}{B}+\varepsilon
\]
for all $s>s_0$.
Since $F(s)<\infty$ and $f>0$ for sufficiently large $s>0$, we have  $F(s)\to 0$ as $s\to \infty$.
Therefore, it holds 
\[
1-\Big(\frac{1}{B}+\varepsilon \Big)\frac{1}{-\log F(s)}
\le 
f'(s)F(s)
\le
1-\Big(\frac{1}{B}-\varepsilon \Big)\frac{1}{-\log F(s)}
\le 1.
\]
This yields that $f'(s)F(s)\rightarrow 1^-$ as $s\to \infty$.
Furthermore, 
$f(s)F(s)$ is non-increasing, namely
\[
\big(f(s)F(s)\big)'=f'(s)F(s)-1\le 0
\]
for all $s>0$.
\end{remark}

\par \medskip
\begin{remark}
The definition of $B_1[f]$
is related to de l'Hospital's form of the definition of $A$.
Let 
\[
A[f](s):=
\frac{\bigl(F(s)\bigr)'}{\left(\frac{1}{f'(s)}\right)'}
=
\frac{f'(s)^2}{f(s)f''(s)}.
\]
Since
\[
(-\log F)'(s)=\frac{1}{f(s)F(s)},
\quad
(-\log F)''(s)=\frac{1-f'(s)F(s)}{\bigl(f(s)F(s)\bigr)^2},
\]
One can observe that
\[
{A\bigl[-\log F\bigr](s)}
=
\frac{(-\log F)'(s)^2}{-\log F(s)(-\log F)''(s)}
=
{B_1[f](s)}.
\]
\end{remark}
\begin{remark}
We point out that in $N \ge 3$ one has for polynomial nonlinearities the ``critical" exponent $p_c := \frac N{N-2}$ (often referred to as {\it Serrin exponent}) which is a borderline exponent: for $ f(u) = u^p$ with $1 < p <  p_c$ the singular solutions behave like the fundamental solution $|x|^{2-N}$, while for $p \ge p_c$ the singular solutions have a weaker growth and become {\it distributional solutions}. For $N = 2$ we have an analogous critical growth given by $f(u) = e^u$. For nonlinearities with growth lower or equal $e^u$ (sub-exponential growth), the singular solutions behave like the fundamental solution $-\log|x|$, while for $f(u) \sim e^{u^q}$ with $q > 1$  (super-exponential growth) the singular solutions have again a weaker growth, and the solutions become distributional.
We remark that in line with these observations, it was proved in \cite{DGP} (see also \cite{GGP}) that for sub-exponential nonlinearities any distributional solution is regular. On the other hand, for super-exponential nonlinearities it is established  that any solution on the punctured ball extends to a  distributional solution, and some  examples of singular distributional solutions are provided.
\end{remark}

\section{Construction of a singular solution}
We first prove Theorem~\ref{theorem:1.2} by constructing a singular solution.
In what follows we consider radially symmetric solutions by focusing on the ordinary differential equation 
\begin{equation}\label{eq:1.1rad}
u''(r)+\frac{N-1}{r}u'(r)+f(u)=0\quad \text{for}\ r>0.
\end{equation}

\subsection{Proof of Theorem~\ref{theorem:1.2}}
Since $v$ defined in \eqref{eq:1.3a}, respectively \eqref{eq:1.3b}, is an explicit solution of $-\Delta v=g(v)$, we define the expected singular function
by
\[
\tilde u(x):=F^{-1}\Bigl[G\bigl(v(x)\bigr)\Bigr].
\]
Then, 
it follows from direct computations as in 
Proposition~3.1 of \cite{FI} 
that $\tilde u(x)$ satisfies
\[
-\Delta {\tilde u}=f(\tilde u)+\frac{|\nabla \tilde{u}|^2}{f(\tilde u)F(\tilde u)}\Bigl[g'(v)G(v)-f'(\tilde u)F(\tilde u)\Bigr].
\]
Since $F(\tilde u)=G(v)$ and $\nabla \tilde u=\frac{f(\tilde u)\nabla v}{g(v)}$, it is equivalent to
\begin{equation}\label{eq:2.1a}
-\Delta {\tilde u}=f(\tilde u)
-f(\tilde u)F(\tilde u)\frac{|\nabla v|^2}{g(v)^2G(v)^2}\Bigl[f'(\tilde u)F(\tilde u)-g'(v)G(v)\Bigr].
\end{equation}
In order to construct a singular solution $u$ satisfying $-\Delta u=f(u)$ near the origin, we look for a solution
in the form
\[
u(x)=\tilde u(x) +\theta(x).
\]
Thanks to \eqref{eq:2.1a}, it is enough to construct a solution $\theta$ satisfying
\begin{equation}\label{eq:2.2a}
-\Delta \theta 
= f(\tilde u + \theta)-f(\tilde u)
+f(\tilde u)F(\tilde u)\frac{|\nabla v|^2}{g(v)^2G(v)^2}\Bigl[f'(\tilde u)F(\tilde u)-g'(v)G(v)\Bigr]
\end{equation}
near the origin.
From now on, we look for a radially symmetric solution to \eqref{eq:2.2a}. 
We make the use of the Emden--Fowler transformation as follows:
\begin{equation}\label{eq:2.0}
\eta(\rho)=\theta(x),\qquad \rho=\log\frac{1}{r^2}+1.
\end{equation}
Then, the equation~\eqref{eq:2.2a} and the functions defined in 
\eqref{eq:1.3a}--\eqref{eq:1.5}
and the condition~\eqref{eq:1.1a}
are transformed into 
\begin{equation}\label{eq:2.1}
\eta''(\rho)
+
\frac{e^{-\rho+1}}{4}\Bigl(f(\phi+\eta)-f(\phi)\Bigr)
+
f(\phi)F(\phi)\frac{\psi'(\rho)^2}{g(\psi)^2G(\psi)^2}
\Bigl[f'(\phi)F(\phi)-g'(\psi)G(\psi)\Bigr]=0,
\end{equation}
and
\begin{align}
&
\label{eq:2.4a}
\psi(\rho):=v(x)=(\rho -1)^{\frac{1}{B'}}, \ {\rm{if}}\ B>1 \ {\rm and} \ \psi(\rho):=v(x)=\log(\rho-1), \ {\rm{if}}\ B=1, \ 
\\
&
\label{eq:2.4bbb}
\phi(\rho):=\tilde u(x)=F^{-1}\Bigl[G\bigl(\psi(\rho)\bigr)\Bigr],
\\
&
\label{eq:2.6b}
\epsilon_1(\rho):=
R_1(x)
=
\bigg|
\frac{1}{B_1[f]\bigl(\phi(\rho)\bigr)}-\frac{1}{B_1[g]\bigl(\psi(\rho)\bigr)}
\bigg|,
\\
&
\label{eq:2.7a}
\epsilon_2(\rho)
:=
R_2(x)
=
\bigg|
\frac{1}{B_2[f]\bigl(\phi(\rho)\bigr)}-\frac{1}{B_2[g]\bigl(\psi(\rho)\bigr)}
\bigg|,
\\
&
\label{eq:5.7}
\lim_{\rho\to \infty}\rho^{\frac{1}{2}}[\epsilon_1(\rho)+\epsilon_2(\rho)]=0.
\end{align}

We decompose the nonlinear term as follows:
\[
\begin{aligned}
\frac{e^{-\rho+1}}{4}\Bigl(f(\phi+\eta)-f(\phi)\Bigr)
=
\frac{e^{-\rho+1}}{4}f'(\phi)\eta
+
\frac{e^{-\rho+1}}{4}\Bigl(f(\phi+\eta)-f(\phi)-f'(\phi)\eta\Bigr).
\end{aligned}
\]
By \eqref{eq:1.4a} and \eqref{eq:2.4bbb},
we have
\begin{equation}\label{eq:FG}
F(\phi)=G(\psi)=\frac{B}{4}\rho e^{-\rho+1}.
\end{equation}
It follows from \eqref{eq:FG} 
that
\[
\begin{aligned}
\frac{e^{-\rho+1}}{4}f'(\phi)\eta
&
=
\frac{e^{-\rho+1}}{4}
\Biggl[
\frac{1}{F(\phi)}
+
\frac{\bigl(-\log F(\phi)\bigr)(f'(\phi)F(\phi)-1)
}{F(\phi)\bigl(-\log F(\phi)\bigr)}
\Biggr]\eta
\\
&
=
\frac{1}{B\rho}\eta+
\frac{e^{-\rho+1}}{4}
\frac{\bigl(-\log F(\phi)\bigr)(f'(\phi)F(\phi)-1)}{F(\phi)\bigl(-\log F(\phi)\bigr)}
\eta.
\end{aligned}
\]
Define
\begin{equation}\label{eq:2.10a}
\left\{
\begin{aligned}
&
I(\rho):=
f(\phi)F(\phi)\frac{\psi'(\rho)^2}{g(\psi)^2G(\phi)^2}
\Bigl[f'(\phi)F(\phi)-g'(\psi)G(\psi)\Bigr],
\\
&
L(\rho):=-\frac{3}{16\rho^2}+\frac{e^{-\rho+1}}{4}
\frac{\bigl(-\log F(\phi)\bigr)(f'(\phi)F(\phi)-1)}{F(\phi)\bigl(-\log F(\phi)\bigr)},
\\
&
N[\eta](\rho):=
\frac{e^{-\rho+1}}{4}\Bigl(f(\phi+\eta)-f(\phi)-f'(\phi)\eta\Bigr),
\end{aligned}
\right.
\end{equation}
then one can rewrite the equation of $\eta$ to the following:
\begin{equation}\label{eq:2.8a}
\eta''(\rho)+\left(\frac{1}{B\rho}+\frac{3}{16\rho^2}\right)\eta
+I(\rho)+L(\rho)\eta+N[\eta](\rho)=0.
\end{equation}
It is easy to check that the functions $(\Phi,\Psi)$ defined by
\[
\Phi(\rho):=\rho^{\frac{1}{4}}\sin\left(\frac{2}{\sqrt{B}}\rho^{\frac{1}{2}}\right),\quad
\Psi(\rho):=\rho^{\frac{1}{4}}\cos\left(\frac{2}{\sqrt{B}}\rho^{\frac{1}{2}}\right)
\]
are the fundamental system of 
the linear equation
\[
\eta''(\rho)+\left(\frac{1}{B\rho}+\frac{3}{16\rho^2}\right)\eta=0.
\]
It is remarkable that the fundamental solutions can be written by the trigonometric functions thanks to
the modification term $\frac{3}{16\rho^2}\eta$; otherwise the Bessel functions of order 1 appears as the fundamental solutions. This simplification is used in \cite{IKNW}.

Then, the variation of parameters shows us that 
the problem~\eqref{eq:2.8a} with $\eta(\rho)\to 0\ (\rho \to \infty)$ is equivalent to
the following integral equation:
\begin{equation}\label{eq:2.9a}
\begin{aligned}
\eta(\rho)
&
=
\sqrt{B}\int_{\rho}^{\infty}\Bigl(\Phi(\rho)\Psi(\tau)-\Psi(\rho)\Phi(\tau)\Bigr)
\Bigl[I(\tau)+L(\tau)\eta(\tau)+N[\eta](\tau)\Bigr]d\tau
\\
&
=\sqrt{B}
\int_{\rho}^{\infty}(\rho \tau)^{\frac{1}{4}}\sin\left(\frac{2}{\sqrt{B}}(\rho^{\frac{1}{2}}-\tau^{\frac{1}{2}})\right)
\Bigl[I(\tau)+L(\tau)\eta(\tau)+N[\eta](\tau)\Bigr]d\tau
\\
&
=:\mathcal{T}[\eta](\rho).
\end{aligned}
\end{equation}
We need decay estimates of $I$ and $L$ to construct a solution of the integral equation~\eqref{eq:2.9a}.
Since 
$F(\phi)=\frac{B}{4}\rho e^{-\rho+1}$ and
$(-\log F(s))(f'(s)F(s)-1)\to \frac{1}{B}$ as $s\to \infty$,
we have
\begin{equation}\label{eq:2.10c}
L(\rho)=O\Big(\frac{1}{\rho^2}\Big).
\end{equation}
We now estimate the order of the inhomogeneous term $I(\rho)$.
\begin{lemma}\label{lemma:2.1}
Let $I$ be the function defined in \eqref{eq:2.10a}. Then there exists $C>0$ such that
\[
\begin{aligned}
|I(\rho)|\le C \frac{f(\phi)F(\phi)}{\rho}\epsilon_1(\rho)
\quad \text{and}
\quad
|I'(\rho)|\le C \frac{f(\phi)F(\phi)}{\rho^2}\bigl(\epsilon_1(\rho)+\epsilon_2(\rho)\bigr).
\end{aligned}
\]
\end{lemma}

\begin{proof}[Proof of Lemma~\ref{lemma:2.1}]
We first recall that 
\begin{equation}\label{eq:3.16ddd}
\frac{\psi'(\rho)}{g(\psi)G(\psi)}=1-\frac{1}{\rho}
\end{equation}
by \eqref{eq:1.3a},  \eqref{eq:1.3b}, and \eqref{eq:2.4a}.
Then, by
$F(\phi)=G(\psi)
=\frac{B}{4}\rho e^{-\rho+1}
$, we have
\[
\begin{aligned}
I(\rho)
&
=
f(\phi)F(\phi)\frac{\psi'(\rho)^2}{g(\psi)^2G(\psi)^2}
\Bigl[f'(\phi)F(\phi)-g'(\psi)G(\psi)\Bigr]
\\
&
=
f(\phi)F(\phi)\left(1-\frac{1}{\rho}\right)^2
\frac{
\Bigl[(\log F(\phi))\bigl(f'(\phi)F(\phi)-1\bigr)+(\log G(\psi))\bigl(1-g'(\psi)G(\psi)\bigr)\Bigr]
}
{
\log G(\psi)
}
\\
&
=
f(\phi)F(\phi)\left(1-\frac{1}{\rho}\right)^2
\frac{\epsilon_1(\rho)}{-\rho+1-\log \frac{1}{\rho}-\log \frac{4}{B}}.
\end{aligned}
\]
It remains to estimate $I'$.
Recall that by $F(\phi)=G(\phi)$
\[
\phi'(\rho)=\frac{f(\phi)}{g(\psi)}\psi'(\rho)=f(\phi)F(\phi)\frac{\psi'(\rho)}{g(\psi)G(\psi)}
=f(\phi)F(\phi)\left(1-\frac{1}{\rho}\right).
\]
Hence, it holds by simple computations that
\[
\begin{aligned}
\frac{d}{d\rho}\Bigl(f(\phi)F(\phi)\Bigr)
&
=
\Bigl(f'(\phi)F(\phi)-1\Bigr)\phi'
=
\Bigl(f'(\phi)F(\phi)-1\Bigr)\left(1-\frac{1}{\rho}\right)f(\phi)F(\phi)
\\
&
=
\frac{\Bigl(f'(\phi)F(\phi)-1\Bigr)(\log F(\phi))}{\log F(\phi)}\left(1-\frac{1}{\rho}\right)f(\phi)F(\phi)
=
f(\phi)F(\phi)
O\left(\frac{1}{\rho}\right),
\end{aligned}
\]
and
\[
\begin{aligned}
&
\frac{d}{d\rho}\Bigl(f'(\phi)F(\phi)-g'(\psi)G(\psi)\Bigr)
\\
&
=
\left(f''(\phi)F(\phi)-\frac{f'(\phi)}{f(\phi)}\right)\phi'(\rho)
-
\left(g''(\psi)G(\psi)-\frac{g'(\psi)}{g(\psi)}\right)\psi'(\rho)
\\
&
=
\left[
\left(\frac{f(\phi)f''(\phi)F(\phi)}{f'(\phi)}-1\right)f'(\phi)F(\phi)
-
\left(\frac{g(\psi)g''(\psi)G(\psi)}{g'(\psi)}-1\right)g'(\psi)G(\psi)
\right]\Big(1-\frac{1}{\rho}\Big)
\\
&
=
\epsilon_2 (\rho)O\left(\frac{1}{\rho^2}\right).
\end{aligned}
\]
Therefore
\[
I'(\rho)
=
\frac{d}{d\rho}\bigg[f(\phi)F(\phi)\Big(1-\frac{1}{\rho}\Big)^2\bigl(f'(\phi)F(\phi)-g'(\psi)G(\psi)\bigr)\bigg]
=
(\varepsilon_1(\rho)+\varepsilon_2(\rho))f(\phi)F(\phi)O\Big(\frac{1}{\rho^2}\Big).
\]
This completes the proof of Lemma~\ref{lemma:2.1}
\end{proof}

\par \bigskip

By Lemma~\ref{lemma:2.1}, we have the following estimate of the inhomogeneous term.
\begin{lemma}\label{lemma:2.2}
Let $I$ be the function defined in \eqref{eq:2.10a}. Then there exists $C_*>0$ such that
\[
\sqrt{B}\int_{\rho}^{\infty}(\rho \tau)^{\frac{1}{4}}\sin\Big(\frac{2}{\sqrt{B}}(\rho^{\frac{1}{2}}-\tau^{\frac{1}{2}})\Big)
I(\tau)d\tau
\le 
C_* 
f\bigl(\phi(\rho)\bigr)
F\bigl(\phi(\rho)\bigr)
\sup_{\tau \ge \rho}
\Bigl(\epsilon_1(\tau)+\epsilon_2(\tau)\Bigr).
\]
\end{lemma}
\begin{proof}[Proof of Lemma~\ref{lemma:2.2}]
It follows from the integration by parts that
\[
\begin{aligned}
&
\int_{\rho}^{\infty}(\rho \tau)^{\frac{1}{4}}\sin\Big(\frac{2}{\sqrt{B}}(\rho^{\frac{1}{2}}-\tau^{\frac{1}{2}})\Big)
I(\tau)d\tau
\\
&
=
\left[
\sqrt{B}\cos \Big(\frac{2}{\sqrt{B}}(\rho^{\frac{1}{2}}-\tau^{\frac{1}{2}})\Big)
\rho^{\frac{1}{4}}\tau^{\frac34}I(\tau)
\right]_{\tau=\rho}^{\tau=\infty}
\\
&
\qquad 
-\sqrt{B}\rho^{\frac{1}{4}}
\int_{\rho}^{\infty}
\cos \Big(\frac{2}{\sqrt{B}}(\rho^{\frac{1}{2}}-\tau^{\frac{1}{2}})\Big)
\Big(\frac{3}{4}\tau^{-\frac{1}{4}}I(\tau)+\tau^{\frac{3}{4}}I'(\tau)\Big)d\tau.
\end{aligned}
\]
This together with Lemma~\ref{lemma:2.1} yields
\[
\begin{aligned}
&
\left|
\sqrt{B}\int_{\rho}^{\infty}(\rho \tau)^{\frac{1}{4}}\sin\left(\frac{2}{\sqrt{B}}(\rho^{\frac{1}{2}}-\tau^{\frac{1}{2}})\right)
I(\tau)d\tau
\right|
\\
&
\le
C \rho I(\rho)
+
C \rho^{\frac{1}{4}}\int_{\rho}^{\infty}\tau^{-\frac{5}{4}}
f\bigl(\phi(\tau)\bigl)
F\bigl(\phi(\tau)\bigl)
\Bigl(\epsilon_1(\tau)+\epsilon_2(\tau)\Bigr)d\tau
\\
&
\le
C_* 
f\bigl(\phi(\rho)\bigr)
F\bigl(\phi(\rho)\bigr)
\sup_{\tau \ge \rho}
\Bigl(\epsilon_1(\tau)+\epsilon_2(\tau)\Bigr).
\end{aligned}
\]
for some constants $C,C_*>0$.
Remark that $f(\phi)F(\phi)$ is a non-increasing function
since 
\[
\frac{d}{d\rho}f(\phi)F(\phi)=[f'(\phi)F(\phi)-1]\phi'(\rho)
=
[f'(\phi)F(\phi)-1]\frac{f(\phi)}{g(\psi)}\psi'
\le 0.
\]
This completes the proof of Lemma~\ref{lemma:2.2}.
\end{proof}
\par \bigskip
We now prepare an estimate of the nonlinear term $N[\eta](\rho)$.
\begin{lemma}\label{lemma:2.3}
Suppose that 
\[
\eta_1(\rho),\eta_2(\rho)
=
O\left(
f\bigl(\phi(\rho)\bigr)
F\bigl(\phi(\rho)\bigr)
\sup_{\tau \ge \rho}
\Bigl(\epsilon_1(\tau)+\epsilon_2(\tau)\Bigr)
\right).
\]
Then, there exist $\rho_0>0$ and $C>0$ such that for all $\rho\ge \rho_0$ it holds
\[
\left|
N[\eta_1](\rho)
-
N[\eta_2](\rho)
\right|
\le 
C 
\frac{1}{\rho}
\sup_{\tau\ge \rho}\bigl(\epsilon_1(\tau)+\epsilon_2(\tau)\bigr)
|\eta_1(\rho)-\eta_2(\rho)|.
\]
\end{lemma}
\begin{proof}[Proof of Lemma~\ref{lemma:2.3}]
Taylor's expansion and $F(\phi(\rho))=\frac{B}{4}\rho e^{-\rho +1}$ show us that
\begin{equation}\label{eq:2.16a}
\begin{aligned}
N[\eta_1](\rho)
-
N[\eta_2](\rho)
&
=
\frac{e^{-\rho+1}}{4}
f(\phi+\eta_1)-f(\phi+\eta_2)-f'(\phi)(\eta_1-\eta_2)
\\
&
=
\frac{e^{-\rho+1}}{4}
\Bigl(f'(\phi+\tilde \eta)-f'(\phi)\Bigr)(\eta_1-\eta_2)
\\
&
=
\frac{e^{-\rho+1}}{4}
f''(\phi+c\tilde \eta)\tilde \eta \bigl(\eta_1-\eta_2\bigr)
\\
&
=
\frac{F(\phi(\rho))}{B\rho}
f''(\phi+c\tilde \eta)\tilde \eta \bigl(\eta_1-\eta_2\bigr),
\end{aligned}
\end{equation}
where $\tilde \eta=a\eta_1+(1-a)\eta_2$ for some $a\in (0,1)$
and $c\in (0,1)$.
It is enough to prove that there exists $\rho_0>0$ such that for all $\rho\ge \rho_0$ it holds
\begin{equation}\label{eq:2.15a}
\frac{1}{2}
\frac{1}{f(\phi)F(\phi)^2}
\le
f''(\phi+c\tilde \eta)
\le 
2\frac{1}{f(\phi)F(\phi)^2}.
\end{equation}
Indeed the conclusion follows from the inequality \eqref{eq:2.15a} combining with \eqref{eq:2.16a} 
and the fact that $\tilde \eta$ has the same order as $\eta_1$ and $\eta_2$.

We start by decomposing $f''$ to
\[
f''=\frac{f f'' F}{f'}\frac{f'F}{fF^2}.
\]
By the assumption \eqref{eq:1.2c}, $\lim_{s\to \infty}f'(s)F(s)=1$, 
and $\lim_{s\to \infty}(-\log F(s))=\infty$, we have
$\lim_{s\to \infty}\frac{f(s)f''(s)F(s)}{f'(s)}=1$.
Moreover, by $\lim_{\rho \to \infty}\phi(\rho)=\infty$ and $\lim_{\rho\to \infty}\tilde \eta(\rho)=0$,
it holds
\[
\phi(\rho)+c\tilde \eta(\rho)\ge \phi - |\tilde \eta|\to \infty
\]
as $\rho \to \infty$. All together, 
there exists $\rho_0>0$ such that for all $\rho\ge \rho_0$ 
it holds
\[
\frac{1}{2}
\frac{1}{f(\phi+c\tilde \eta)F(\phi+c\tilde \eta)^2}
\le 
f''(\phi+c\tilde \eta) 
\le
2
\frac{1}{f(\phi+c\tilde \eta)F(\phi+c\tilde \eta)^2}.
\]
Therefore, it remains to prove \eqref{eq:2.15a}
that
\begin{equation}\label{eq:2.17a}
\lim_{\rho\to \infty}\frac{f(\phi+c\tilde \eta)F(\phi+c\tilde \eta)}{f(\phi)F(\phi)}=1
\end{equation}
and
\begin{equation}\label{eq:2.18a}
\lim_{\rho\to \infty}\frac{F(\phi+c\tilde \eta)}{F(\phi)}=1.
\end{equation}
By the Taylor expansion, there exists $c'\in (0,1)$
such that
\[
f(\phi+c\tilde \eta)F(\phi+c\tilde \eta)=
f(\phi)F(\phi)+\Bigl[f'(\phi+c'\tilde \eta)F(\phi+c'\tilde \eta)-1\Bigr]c\tilde \eta.
\]
Hence, 
\[
\frac{
f(\phi+c\tilde \eta)F(\phi+c\tilde \eta)
}
{
f(\phi)F(\phi)
}
=
1+\Bigl[f'(\phi+c'\tilde \eta)F(\phi+c'\tilde \eta)-1\Bigr]c
\frac{\tilde \eta}{f(\phi)F(\phi)}.
\]
By the assumption on $\eta_1$ and $\eta_2$, we see 
\begin{equation}\label{eq:2.19a}
\left|
\frac{\tilde \eta}{f(\phi)F(\phi)}
\right|
\le
C
\sup_{\tau \ge \rho}
\Bigl(\epsilon_1(\tau)+\epsilon_2(\tau)\Bigr)
\to 0
\end{equation}
as $\rho \to \infty$. 
Moreover, $f'(s)F(s)\to 1$ as $s\to \infty$.
This proves \eqref{eq:2.17a}.

We turn to prove \eqref{eq:2.18a}. 
By the Taylor expansion
and the fact that $f$ is increasing, there exists $c_1,c_2\in (0,1)$ such that
\[
\begin{aligned}
&
F(\phi+|\tilde \eta|)
=
F(\phi)-\frac{1}{f(\phi+c_1|\tilde \eta|)}|\tilde \eta|
\ge 
F(\phi)-\frac{1}{f(\phi)}|\tilde \eta|,
\\
&
F(\phi-|\tilde \eta|)=F(\phi)+\frac{1}{f(\phi-c_2|\tilde \eta|)}|\tilde \eta|
\le F(\phi)+\frac{|\tilde \eta|}{f(\phi-|\tilde \eta|)}.
\end{aligned}
\]
Combining \eqref{eq:2.17a} and \eqref{eq:2.19a}, we see 
$\frac{|\tilde \eta|}{f(\phi-|\tilde \eta|)F(\phi-|\tilde \eta|)}\to 0$ as $\rho \to \infty$. 
Since
$F$ is decreasing, 
we obtain
\[
1-\frac{|\tilde \eta|}{f(\phi)F(\phi)}
\le 
\frac{F(\phi+|\tilde \eta|)}{F(\phi)}
\le 
\frac{F(\phi)}{F(\phi)}=1.
\]
This yields 
\[
\lim_{\rho \to \infty}\frac{F(\phi+|\tilde \eta|)}{F(\phi)}
=
\lim_{\rho \to \infty}\frac{F(\phi-|\tilde \eta|)}{F(\phi)}
=1.
\]
Applying these limits to $\frac{F(\phi+|\tilde \eta|)}{F(\phi)}\le \frac{F(\phi+c\tilde \eta)}{F(\phi)}\leq
\frac{F(\phi-|\tilde \eta|)}{F(\phi)}$, we obtain \eqref{eq:2.18a}.
This completes the proof of Lemma~\ref{lemma:2.3}.
\end{proof}

\par \bigskip
We are now in a position to 
construct a fixed point of the operator $\mathcal{T}$ defined in \eqref{eq:2.9a}.
Let $C_*>0$ be the constant taken in Lemma~\ref{lemma:2.2}.
Fix $\rho_{0}>0$.
By Lemma~\ref{lemma:2.2}, we defined the space $X_{\rho_0}$ by
\[
\begin{aligned}
X_{\rho_0}
&
:=
\Bigl\{
\eta\in C\Bigl([\rho_0,\infty)\Bigr):\|\eta\|
\le 
2C_*
\Bigr\},
\\
\|\eta\|
&
:=
\sup_{\rho\ge \rho_0}
\frac{
|\eta(\rho)|
}
{
f\bigl(\phi(\rho)\bigr)F\bigl(\phi(\rho)\bigr)\sup_{\tau \ge \rho}\Bigl(\epsilon_1(\tau)+\epsilon_2(\tau)\Bigr)
}.
\end{aligned}
\]

We prove that there exists sufficiently large $\rho_0>0$ such that $\mathcal{T}$ is a contraction map on $X_{\rho_0}$.
Let $\eta \in X_{\rho_0}$, where $\rho_0$ will be chosen later. 
By \eqref{eq:2.10c} we have
\[
\begin{aligned}
\left|
\int_{\rho}^{\infty}(\rho \tau)^{\frac{1}{4}}\sin\Big(\frac{2}{\sqrt{B}}(\rho^{\frac{1}{2}}-\tau^{\frac{1}{2}})\Big)
L(\tau)\eta(\tau)d\tau
\right|
&
\le
\int_{\rho}^{\infty}(\rho \tau)^{\frac{1}{4}}
\frac{C}{\tau^2}|\eta(\tau)|d\tau
\\
&
\le 
C \rho^{-\frac{1}{2}}
f(\phi(\rho))F(\phi(\rho))\sup_{\tau\ge \rho}(\epsilon_1(\tau)+\epsilon_2(\tau)).
\end{aligned}
\]
Taking $\rho_0>0$ such that $C \rho_0^{-\frac{1}{2}}\le \frac{C_*}{2}$,
we obtain 
\[
\left\|
\int_{\rho}^{\infty}(\rho \tau)^{\frac{1}{4}}\sin\left(\frac{2}{\sqrt{B}}(\rho^{\frac{1}{2}}-\tau^{\frac{1}{2}})\right)
L(\tau)\eta(\tau)d\tau
\right\|\le \frac{C_*}{2}.
\]
On the other hand,
Lemma~\ref{lemma:2.3} with $\eta_1=\eta$ and $\eta_2=0$
yields that
\[
\begin{aligned}
&
\left|
\int_{\rho}^{\infty}(\rho \tau)^{\frac{1}{4}}\sin\Big(\frac{2}{\sqrt{B}}(\rho^{\frac{1}{2}}-\tau^{\frac{1}{2}})\Big)
N[\eta](\tau)d\tau
\right|
\le 
C \int_{\rho}^{\infty}(\rho \tau)^{\frac{1}{4}}
\frac{1}{\tau}
\sup_{t\ge \tau}\bigl(\epsilon_1(t)+\epsilon_2(t)\bigr)
\eta(\tau)
d\tau.
\end{aligned}
\]
Applying the assumption \eqref{eq:1.1a}, equivalently \eqref{eq:5.7},
we obtain
\[
\begin{aligned}
&
\int_{\rho}^{\infty}(\rho \tau)^{\frac{1}{4}}
\frac{1}{\tau}
\sup_{t\ge \tau}\bigl(\epsilon_1(t)+\epsilon_2(t)\bigr)
\eta(\tau)
d\tau
\\
&
\le 
C 
\Bigl(\sup_{\tau \ge \rho}\tau^{\frac{1}{2}}[\epsilon_1(\tau)+\epsilon_2(\tau)]\Bigr)
\Bigl(
f(\phi(\rho))F(\phi(\rho))\sup_{\tau\ge \rho}[\epsilon_1(\tau)+\epsilon_2(\tau)]
\Bigr)
\int_{\rho}^{\infty}(\rho\tau)^{\frac{1}{4}}\frac{1}{\tau^{\frac{1}{2}}}\frac{1}{\tau}d\tau
\\
&
\le 
C 
\Bigl(\sup_{\tau \ge \rho}\tau^{\frac{1}{2}}[\epsilon_1(\tau)+\epsilon_2(\tau)]\Bigr)
\Bigl(
f(\phi(\rho))F(\phi(\rho))\sup_{\tau\ge \rho}[\epsilon_1(\tau)+\epsilon_2(\tau)]
\Bigr).
\end{aligned}
\]
This together with the assumption \eqref{eq:1.1a}, one can take $\rho_0$ such that
\[
\left\|
\int_{\rho}^{\infty}(\rho \tau)^{\frac{1}{4}}\sin\Big(\frac{2}{\sqrt{B}}(\rho^{\frac{1}{2}}-\tau^{\frac{1}{2}})\Big)
N[\eta](\tau)d\tau
\right\|
\le 
\frac{C_*}{2}
\]
for all $\rho\ge \rho_0$.
This proves that $\mathcal{T}[\eta]\in X_{\rho_0}$. 

Next we prove that $\mathcal{T}$ is a contraction mapping.
For $\eta_1,\eta_2\in X_{\rho_0}$,
we have
\begin{equation}\label{eq:2.18b}
\begin{aligned}
|\mathcal{T}[\eta_1](\rho)-\mathcal{T}[\eta_2](\rho)|
\le 
\sqrt{B}\int_{\rho}^{\infty}(\rho \tau)^{\frac{1}{4}}
\biggl[|L(\tau)||\eta_1-\eta_2|+\Bigl|N[\eta_1](\tau)-N[\eta_2](\tau)\Bigr|\biggr]d\tau.
\end{aligned}
\end{equation}
By \eqref{eq:2.10c}, it holds
\begin{equation}\label{eq:2.21a}
\begin{aligned}
&
\sqrt{B}\int_{\rho}^{\infty}(\rho \tau)^{\frac{1}{4}}
|L(\tau)||\eta_1-\eta_2|d\tau
\\
&
\le
C
\int_{\rho}^{\infty}\rho^{\frac{1}{4}}\tau^{\frac{1}{4}-2}d\tau
\Bigl(
f(\phi(\rho))F(\phi(\rho))\sup_{\tau\ge \rho}(\epsilon_1(\tau)+\epsilon_2(\tau))
\Bigr)
\|\eta_1-\eta_2\|
\\
&
\le
\frac{1}{3}
\Bigl(
f(\phi(\rho))F(\phi(\rho))\sup_{\tau\ge \rho}(\epsilon_1(\tau)+\epsilon_2(\tau))
\Bigr)
\|\eta_1-\eta_2\|
\end{aligned}
\end{equation}
for large $\rho>0$.
Since $fF$ is non-increasing by Remark~\ref{remark:1.3}, it follows from Lemma~\ref{lemma:2.3} and the assumption~\eqref{eq:1.1a} that
\begin{equation}\label{eq:2.20a}
\begin{aligned}
&
\sqrt{B}\int_{\rho}^{\infty}(\rho \tau)^{\frac{1}{4}}
\Bigl|N[\eta_1](\tau)-N[\eta_2](\tau)\Bigr|d\tau
\\
&
\le
C
\int_{\rho}^{\infty}(\rho \tau)^{\frac{1}{4}}
\frac{1}{\tau}
\sup_{t\ge \tau}\bigl(\epsilon_1(t)+\epsilon_2(t)\bigr)
|\eta_1(\tau)-\eta_2(\tau)|d\tau
\\
&
=
C
\rho^{\frac{1}{4}}\int_{\rho}^{\infty}\tau^{\frac{1}{4}-\frac{1}{2}-1}
\Big(\tau^{\frac{1}{2}}\sup_{t\ge \tau}\bigl(\epsilon_1(t)+\epsilon_2(t)\bigr)\Big)
\frac{f(\tau)F(\tau)\sup_{t\ge \tau}\bigl(\epsilon_1(t)+\epsilon_2(t)\bigr)}{f(\tau)F(\tau)\sup_{t\ge \tau}\bigl(\epsilon_1(t)+\epsilon_2(t)\bigr)}
|\eta_1(\tau)-\eta_2(\tau)|d\tau
\\
&
\le
C 
\sup_{\tau \ge \rho}\Bigl(\tau^{\frac{1}{2}}\sup_{t\ge \tau}\bigl(\epsilon_1(t)+\epsilon_2(t)\bigr)\Bigr)
\Bigl(
f(\phi(\rho))F(\phi(\rho))\sup_{t\ge \rho}(\epsilon_1(t)+\epsilon_2(t))
\Bigr)
\|\eta_1-\eta_2\|
\Big(\rho^{\frac{1}{4}}\int_{\rho}^{\infty}\tau^{-\frac{5}{4}}d\tau\Big)
\\
&
\le
\frac{1}{3}
\Bigl(
f(\phi(\rho))F(\phi(\rho))\sup_{t\ge \rho}(\epsilon_1(t)+\epsilon_2(t))
\Bigr)
\|\eta_1-\eta_2\|
\end{aligned}
\end{equation}
for sufficiently large $\rho>0$.
Combining all estimates in \eqref{eq:2.18b}--\eqref{eq:2.20a}, we obtain
\[
\bigl\|\mathcal{T}[\eta_1]-\mathcal{T}[\eta_2]\bigr\|
\le 
\frac{2}{3}\|\eta_1-\eta_2\|
\]
for sufficiently large $\rho$.
Then the contraction mapping argument gives that there exists $\rho_0>0$ and $\eta\in X_{\rho_0}$
such that
$\eta=\mathcal{T}[\eta]$ on $(\rho_0,\infty)$.
Consequently, 
this completes the proof of Theorem~\ref{theorem:1.2}.

\subsection{Proof of Theorem~\ref{theorem:1.1}}
We extend the solution for $r>r_0$ by the shooting methods so that the 0-Dirichlet boundary condition is satisfied.
\begin{proposition}
\label{proposition:5.1}
Assume that $f$ satisfies all assumptions as in Theorem~\ref{theorem:1.1}.
Let $u$ be the solution constructed in Theorem~\ref{theorem:1.2} and $w$ be a solution of
\begin{equation}
\label{eq:5.1}
\left\{ 
\begin{array}{ll}
-w'' - \frac 1r w' = f(w)  \hbox{ for } r > r_0,
\vspace{0.2cm} \\
w(r_0) = u(r_0),
\vspace{0.2cm} \\
w'(r_0) = u'(r_0).
\end{array}
\right.
\end{equation}
Then there exists $R > r_0$ such that $w(R) = 0$.
\end{proposition}
Theorem~\ref{theorem:1.2} and Proposition~\ref{proposition:5.1}
give us that 
the function $u_{\infty}$ defined by 
\begin{equation}\label{eq:2.2aaaaa}
u_{\infty}(x):=
\left\{
\begin{aligned}
&
u(x)
&&\text{in}\ \  B_{r_0}\setminus \{0\},
\\
&
w(|x|)
&&\text{in}\ \  B_R\setminus B_{r_0},
\end{aligned}
\right.
\end{equation}
satisfies 
\eqref{eq:1.1aaa} in the classical sense. 

\par \bigskip

The following lemma is used in the proof of Proposition~\ref{proposition:5.1}.
\begin{lemma}\label{lemma:5.1}
Let $u$ be the solution constructed in Theorem~\ref{theorem:1.2} and $v,g$ be functions defined in \eqref{eq:1.3a}  
or \eqref{eq:1.3b}.
Then it holds
\[
\lim_{r\to 0} \frac{u'(r)}{f\bigl(u(r)\bigr)}\left\slash \frac{v'(r)}{g\bigl(v(r)\bigr)}\right.=1.
\]
\end{lemma}
\begin{proof}
By the change of variables in \eqref{eq:2.0}, it holds
\[
u(r)=\phi(\rho)+\eta(\rho),\  v(r)=\psi(\rho)
\]
and hence
\[
\frac{u'(r)}{f(u(r))}=\frac{\frac{d\rho}{dr}(\phi'+\eta')}{f(\phi+\eta)}=-\frac{2}{r}
\frac{\phi'+\eta'}{f\bigl(\phi+\eta\bigr)}
.
\]
By Taylor's expansion of $f$, there exists $0<c<1$ such that
\[
\frac{\phi'+\eta'}{f\bigl(\phi+\eta\bigr)}
=
\frac{\phi'+\eta'}{f(\phi)+f'(\phi+c\eta)\eta}
=
\frac{\phi'}{f(\phi)}
\times \frac{1+\frac{\eta'}{\phi'}}{1+\frac{f'(\phi+c\eta)}{f(\phi)}\eta}.
\]
Since
$
\frac{\phi'}{f(\phi)}
=
\frac{\psi'}{g(\psi)}
$
and
$\frac{v'(r)}{g(v(r))}=-\frac{2}{r}\frac{\psi'(\rho)}{g(\psi(\rho))}$,
it remains to prove 
that
\begin{equation}
\label{eq:5.3}
\lim_{\rho \to \infty}\frac{f'\bigl(\phi(\rho)+c\eta(\rho)\bigr)}{f(\phi(\rho))}\eta(\rho)=0
\quad 
\text{and}
\quad
\lim_{\rho\to \infty}\frac{\eta'(\rho)}{\phi'(\rho)}=0.
\end{equation}
By the completely same argument 
 as in the derivation of \eqref{eq:2.18a}, we obtain
$\lim_{\rho\to \infty}\frac{F(\phi+c\eta)}{F(\phi)}=1$.
This and Theorem~\ref{theorem:1.2}, and $\lim_{s\to \infty}f'(s)F(s)= 1$
 yield that
\[
\lim_{\rho \to \infty}\frac{f'\bigl(\phi(\rho)+c\eta(\rho)\bigr)}{f(\phi(\rho))}\eta(\rho)
=
\lim_{\rho \to \infty}\frac{F(\phi)}{F(\phi+c\eta)}\times \frac{\eta}{f(\phi)F(\phi)}\times f'\bigl(\phi+c\eta\bigr)F\bigl(\phi+c\eta\bigr)
=
1\times 0 \times 1=0.
\]
We now estimate $\phi'$ and $\eta'$ to show the second limit in \eqref{eq:5.3}.
Since $F(\phi)=G(\psi)$
 and \eqref{eq:3.16ddd}, we have
\[
\phi'
=
f(\phi)\frac{\psi'}{g(\psi)}
=
f(\phi)F(\phi)\frac{\psi'}{g(\psi)G(\psi)}
=
f(\phi)F(\phi)\left(1-\frac{1}{\rho}\right).
\]
On the other hand, since $\eta$ satisfies the integral equation~\eqref{eq:2.9a}, it holds
\[
\eta'(\rho)=\frac{\eta(\rho)}{4\rho}
+
\sqrt{B}
\int_\rho^{\infty}(\rho\tau)^{\frac{1}{4}}
\cos\left(\frac{2}{\sqrt{B}}(\rho^{\frac{1}{2}}-\tau^{\frac{1}{2}})\right)
\frac{1}{\sqrt{B}}\rho^{-\frac{1}{2}}\Bigl[I(\tau)+L(\tau)\eta(\tau)+N[\eta](\tau)\Bigr]d\tau.
\]
We estimate the terms of $I,L,N$, separately.
We start from the estimate of $I$.
Integrating by parts, we have from Lemma~\ref{lemma:2.1} that
\begin{equation}\label{eq:5.4}
\begin{aligned}
&
\left|
\int_\rho^{\infty}
(\rho\tau)^{\frac{1}{4}}
\cos\left(\frac{2}{\sqrt{B}}(\rho^{\frac{1}{2}}-\tau^{\frac{1}{2}})\right)
\rho^{-\frac{1}{2}}I(\tau)
d\tau
\right|
\\
&
=
\sqrt{B}
\left|
\int_{\rho}^{\infty}
\rho^{-\frac{1}{4}}
\sin\left(\frac{2}{\sqrt{B}}(\rho^{\frac{1}{2}}-\tau^{\frac{1}{2}})\right)
{
\left(\tau^{\frac{3}{4}}I(\tau)\right)'
}
d\tau
\right|
\\
&
\le
\sqrt{B}\rho^{-\frac{1}{4}}\int_{\rho}^{\infty}
\left|
{
\frac{3}{4}\tau^{-\frac{1}{4}}I(\tau)+\tau^{\frac{3}{4}}I'(\tau)
}
\right|
d\tau
\\
&
\le 
C
\rho^{-\frac{1}{4}}\int_{\rho}^{\infty}
\tau^{-\frac{5}{4}}f(\phi(\tau))F(\phi(\tau))
\sup_{t\ge \tau}
\Bigl(\varepsilon_1(t)+\varepsilon_2(t)\Bigr)
d\tau
\\
&
\le 
C
\rho^{-\frac{1}{2}}f(\phi(\rho))F(\phi(\rho))
\sup_{\tau\ge \rho}
\Bigl(\varepsilon_1(\tau)+\varepsilon_2(\tau)\Bigr)
\end{aligned}
\end{equation}
for some $C>0$.
The term of $L$ can be estimated from \eqref{eq:2.10c} and Theorem~\ref{theorem:1.2} 
as follows:
\begin{equation}\label{eq:5.5}
\begin{aligned}
&
\left|
\int_\rho^{\infty}
(\rho\tau)^{\frac{1}{4}}
\cos\left(\frac{2}{\sqrt{B}}(\rho^{\frac{1}{2}}-\tau^{\frac{1}{2}})\right)
\rho^{-\frac{1}{2}}L(\tau)\eta(\tau)
d\tau
\right|
\\
&
\le
C\rho^{-\frac{1}{4}}
\sup_{\tau\ge \rho}|\eta(\tau)|\int_{\rho}^{\infty}\tau^{-\frac{7}{4}}d\tau
\\
&
\le
C \frac{1}{\rho}f(\phi(\rho))F(\phi(\rho))
\sup_{\tau\ge \rho}
\Bigl(\varepsilon_1(\tau)+\varepsilon_2(\tau)\Bigr)
\end{aligned}
\end{equation}
for some $C>0$.
Finally, by Lemma~\ref{lemma:2.3} and \eqref{eq:5.7} it holds
\begin{equation}\label{eq:5.6}
\begin{aligned}
&
\left|
\int_\rho^{\infty}
(\rho\tau)^{\frac{1}{4}}
\cos\left(\frac{2}{\sqrt{B}}(\rho^{\frac{1}{2}}-\tau^{\frac{1}{2}})\right)
\rho^{-\frac{1}{2}}N[\eta](\tau)
d\tau
\right|
\\
&
\le
C\rho^{-\frac{1}{4}}\int_{\rho}^{\infty}\tau^{\frac{1}{4}-1}
\sup_{t\ge \tau}\Bigl(\epsilon_1(t)+\epsilon_2(t)\Bigr)|\eta(\tau)|d\tau
\\
&
\le 
C\rho^{-\frac{1}{4}}\sup_{\tau \ge \rho}|\eta(\tau)|\int_{\rho}^{\infty}\tau^{-\frac{5}{4}}
\left[
\tau^{\frac{1}{2}}\sup_{t\ge \tau}\Bigl(\epsilon_1(t)+\epsilon_2(t)\Bigr)
\right]
d\tau
\\
&
\le C\rho^{-\frac{1}{2}}f(\phi(\rho))F(\phi(\rho))
\sup_{\tau\ge \rho}\Bigl(\epsilon_1(\tau)+\epsilon_2(\tau)\Bigr).
\end{aligned}
\end{equation}
Combining \eqref{eq:5.4}--\eqref{eq:5.6}, we obtain
\[
\left|
\frac{\eta'(\rho)}{\phi'(\rho)}
\right|
\le C \rho^{-\frac{1}{2}}
\sup_{\tau\ge \rho}\Bigl(\epsilon_1(\tau)+\epsilon_2(\tau)\Bigr)
\to 0\ (\rho\to \infty).
\]
This shows \eqref{eq:5.3} and hence Lemma~\ref{lemma:5.1} is proved. 
\end{proof}

\par \medskip \noindent
\begin{proof}[Proof of Proposition~\ref{proposition:5.1}]
Assume that there is no such $R > r_0$. Then $w(r) > 0$ for all $r > r_0$.
We first show that 
$w'(r) < 0$ for all $r > r_0$.
By Lemma~\ref{lemma:5.1} and the fact that $v'(r)<0$ for small $r>0$, one can choose small $r_0$
such that $w'(r_0)=u'(r_0)<0$. 
If $w'(r) \ge 0$ for some $r > r_0$, 
there exists $r_1 > r_0$ such
that $w'(r_1) = 0$
and
$w'(r)<0$ for all $r_0\le r < r_1$. 
Then $w''(r_1)=\lim_{\epsilon \to 0+}\frac{w'(r_1)-w'(r_1-\epsilon)}{\epsilon} \ge 0$,
but by the equation~\eqref{eq:5.1} it holds $w''(r_1) = - f(w(r_1))<0$.
Hence $w'(r) < 0$ for all $r > r_0$.

We then prove that there exists $a>0$ such that
$\lim_{r\to \infty}rw'(r)=-a$.
To this end, we show that there exists $C\in \mathbb{R}$ satisfying
\begin{equation}\label{eq:5.2a}
F(w(r))-\frac{r^2}{4}
\ge 
C
\end{equation}
for all $r>r_0$.
Since $w$ is decreasing, we have $w'<0$, and using $f' >0$
\begin{equation}
\notag
\begin{aligned}
-rw'(r)
&
=\int_{r_0}^r sf(w(s))ds -r_0w'(r_0)
\\
&
\ge
f(w(r))\left(\frac{r^2-r_0^2}{2}\right)-r_0w'(r_0)
\\
&
\ge 
f(w(r))\frac{r^2}{2}-f(w(r_0))\frac{r_0^2}{2}-r_0w'(r_0)
\\
&
=
f(w(r))\frac{r^2}{2}-f(u(r_0))\frac{r_0^2}{2}-r_0u'(r_0).
\end{aligned}
\end{equation}
It follows from Lemma~\ref{lemma:5.1}, using $f>0$
\[
\begin{aligned}
-f(u(r_0))\frac{r_0^2}{2}-r_0u'(r_0)
&
=
r_0^2f(u(r_0))\left(-\frac{u'(r_0)}{f(u(r_0))}\frac{1}{r_0}-\frac{1}{2}\right)
\\
&
\ge
r_0^2f(u(r_0))\left(-\frac{1}{2}\frac{v'(r_0)}{g(v(r_0))}\frac{1}{r_0}-\frac{1}{2}\right)
\\
&
=
r_0^2f(u(r_0))\left(\frac{B}{4}\log \frac{1}{r_0^2}-\frac{1}{2}\right)
\ge 0
\end{aligned}
\]
for sufficiently small $r_0>0$,
where $v,g$ are functions defined in \eqref{eq:1.3a}  or \eqref{eq:1.3b}.
This yields that
$-rw'(r)\ge f(w(r))\frac{r^2}{2}$, equivalently, 
$\frac{d}{dr}\Bigl(F(w(r))\Bigr)\ge \frac{r}{2}$.
Hence it holds
\[
F(w(r))-\frac{r^2}{4}
\ge 
F(w(r_0))-\frac{r_0^2}{4}
=
F(u(r_0))-\frac{r_0^2}{4}.
\]
Setting $C=F(u(r_0))-\frac{r_0^2}{4}$, we obtain
 \eqref{eq:5.2a}.

If $F$ is bounded on $(0,\infty)$, namely $\lim_{r\to 0+}F(r)<\infty$, then \eqref{eq:5.2a} immediately yields a contradiction. 
In what follows we assume $\lim_{r\to 0+}F(r)=\infty$.
By \eqref{eq:5.2a}, using $f'>0$, we have 
\[w(r)\le F^{-1}\Big(\frac{r^2}{4}+C\Big)\]
and hence
\[
\int_{r_0}^r sf(w(s))ds
\le
\int_{r_0}^r sf\Big(F^{-1}\Big(\frac{s^2}{4}+C\Big)\Big)ds
=
\Big[-2F^{-1}\Big(\frac{s^2}{4}+C\Big)\Big]_{s=r_0}^{s=r}
\le
2F^{-1}\Big(\frac{r_0^2}{4}+C\Big)
\]
for all $r>r_0$.
Hence 
the integral $\int_{r_0}^{\infty} sf(w(s))ds$ converges. 
Combining this and $w'(r_0)<0$, we obtain
\[
rw'(r)
=
r_0w'(r_0)-\int_{r_0}^r sf(w(s))ds 
\to
r_0w'(r_0)-\int_{r_0}^{\infty} sf(w(s))ds 
<0.
\]
Taking $a=-r_0w'(r_0)+\int_{r_0}^{\infty} sf(w(s))ds>0$, we obtain
$\lim_{r\to \infty}rw'(r)=-a<0$. 
This gives us a contradiction. Indeed,
\[
w(r)=\int_{r_0}^r w'(s)ds+w(r_0)
\le  \int_{r_0}^r -\frac{a}{2s}ds +w(r_0)
=
-\frac{a}{2}\log \frac{r}{r_0}+w(r_0)\to -\infty\ (r\to \infty).
\]
Hence there must exist $R > r_0$ such that $w(R) = 0$. 
\end{proof}

We conclude the proof of Theorem~\ref{theorem:1.1} by proving the following:
\begin{proposition}\label{proposition:2.3}
The function $u_{\infty}$ defined in \eqref{eq:2.2aaaaa} satisfies 
$-\Delta u_{\infty}=f(u_{\infty})$ in $B_R$
in the sense of distributions. 
\end{proposition}
We are also able to estimate the size of the singularity of the solution at the origin.
\begin{lemma}\label{LemfF}
	For any $\varepsilon>0$ there exist $\rho_\varepsilon>0$ and $C>0$ such that for any $\rho>\rho_\varepsilon$
	$$
	\frac{C}{\rho^{\frac 1B+\varepsilon}}\leq f(\phi(\rho)) F(\phi(\rho))\leq \frac{C}{\rho^{\frac 1B-\varepsilon}}
	$$
	where $\phi(\rho)=F^{-1}(G(\psi(\rho))=F^{-1}(\frac{B}4\rho{\rm e}^{-\rho+1})$.
\end{lemma}
\begin{proof}[Proof of Lemma~\ref{LemfF}]
We recall that $F(\phi(\rho))=\frac B4 \rho {\rm e}^{-\rho+1}$  and so we have
$$
\begin{aligned}
\frac{d}{d\rho}\left(f(\phi(\rho))F(\phi(\rho))\right)&=\left(f'(\phi(\rho)) F(\phi(\rho))-1\right) \phi'(\rho)\\
&=\frac{\Big(f'(\phi(\rho)) F(\phi(\rho))-1\Big)\Bigl(-\log {F\bigl(\phi(\rho)\bigr)}\Bigr)}{-\log {F\bigl(\phi(\rho)\bigr)}}\Big(1-\frac 1\rho\Big) \frac B4 \rho {\rm e}^{-\rho+1}  f(\phi(\rho))\\
&=\frac{1}{B_1[f]\bigl(\phi(\rho)\bigr)}\frac{1-\frac 1\rho}{\rho-1-\log\rho-\log\frac B4}f(\phi(\rho))F(\phi(\rho)).
\end{aligned}
$$
Since    $\frac{1}{B_1[f](\phi(\rho))}\to \frac 1B$ as $\rho\to \infty$, we have that
for any $\varepsilon >0$  there exists $\rho_\varepsilon>0$ such that for any $\rho>\rho_\varepsilon$
$$
-\Big(\frac 1B+\varepsilon\Big)\frac1\rho\leq \frac{\frac{d}{d\rho}\left(f(\phi(\rho))F(\phi(\rho))\right)}{f(\phi(\rho))F(\phi(\rho))}\leq-\Big(\frac 1B-\varepsilon\Big)\frac1\rho
$$
and integrating over $(\rho,\infty)$ we have
$$
\frac{C}{\rho^{\frac 1B+\varepsilon}}\leq
f(\phi(\rho))F(\phi(\rho))\leq\frac{C}{\rho^{\frac 1B-\varepsilon}}
$$
for some $C>0$. 
This completes the proof.
\end{proof}
Lemma~\ref{LemfF} implies the following asymptotic behavior of the singular solution $u_{\infty}$.
\begin{lemma}\label{est}
	For any $\sigma >0$ there exists $r_\sigma$ such that
\begin{equation}\label{eq:f}
\frac{C}{|x|^2 \left(1-2\log|x|\right)^{1+\frac{1}{B}+\sigma}}
\le 
f\bigl(u_{\infty}(x)\bigr)\le
\frac{C}{|x|^2 \left(1-2\log|x|\right)^{1+\frac{1}{B}-\sigma}},
\end{equation}
	\begin{equation}\label{Stima1}
	0\leq  u_\infty(x)\leq (-2\log|x|)^{1-\frac 1B+\sigma},
	\end{equation}
	and
		\begin{equation}\label{Stima2}
	|\nabla u_\infty(x)|\leq \frac{1}{|x|(-2\log|x|)^{\frac 1B-\sigma}}
\end{equation}
	for all $x$ with
	$|x|\leq r_\sigma$.
\end{lemma}
\begin{remark}
As a consequence of Lemma~\ref{est}, we have $u_{\infty}\in L^p$ for all $p\in [1,\infty)$ and
$f(u_{\infty})\in L^1(B_R)$ and $f(u_{\infty})\notin L^p(B_R)$
for any $p\in (1,\infty]$.
Furthermore,
if $1<B<2$ then there exists small $r_0>0$ such that the singular solution 
$u_{\infty}$ belongs to the energy class $H^1(B_{r_0})$.
\end{remark}
\begin{proof}[Proof of Lemma~\ref{est}]
Since $u_\infty(x)\sim  \tilde u(x)$ as $x\to0$, in order to establish  the inequality \eqref{Stima1},  it suffices to prove that for any  $\sigma>0$ there exists $r_\sigma >0$ such that
\begin{equation}\label{tilde1}
0\leq \tilde u(x)\leq  (-2\log|x|)^{1-\frac 1B+\sigma}
\end{equation}
for any  $x\in B_R$ with $|x|\leq r_\sigma.$  
 Now, thanks to the transformation $\rho=\log\frac 1{r^2}+1,$ with $r=|x|,$ the inequality \eqref{tilde1} is a consequence of
\begin{equation}\label{tilde11}
\lim_{\rho\to \infty}\frac{\phi(\rho)}{(\rho-1)^{1-\frac 1B+\sigma}}=0.
\end{equation}
Since $\phi(\rho)=F^{-1}\left(\frac B4 \rho {\rm e}^{-\rho+1}\right)$, we have by de l'Hospital's rule
\begin{equation*}
\begin{aligned}
\lim_{\rho\to \infty}\frac{\phi(\rho)}{(\rho-1)^{1-\frac 1B+\sigma}}&
=\lim_{\rho\to \infty}\frac{F^{-1}\left(\frac B4 \rho {\rm e}^{-\rho+1}\right)}{(\rho-1)^{1-\frac 1B+\sigma}}\\
&=\lim_{\rho\to \infty}\frac{f\left(F^{-1}\left(\frac B4 \rho {\rm e}^{-\rho+1}\right)\right) \frac B4 {\rm e}^{-\rho+1} (\rho-1)}{ (1-\frac 1B+\sigma)(\rho-1)^{-\frac 1B+\sigma}}\\
&=\lim_{\rho\to \infty}\frac{f\left(\phi(\rho)\right) F(\phi(\rho))\left(1-\frac 1\rho\right)}  { (1-\frac 1B+\sigma)(\rho-1)^{-\frac 1B+\sigma}}.
\end{aligned}
\end{equation*}	
Then, thanks to Lemma \ref{LemfF}, for all $\varepsilon>0$ with  $\rho>\rho_\varepsilon$, we have
$$
0\leq f(\phi(\rho)) F(\phi(\rho))\leq \frac{C}{\rho^{\frac 1B-\varepsilon}}
$$
for some $C>0$. Thus, choosing $\varepsilon <\sigma$,  for  $\rho>\rho_\varepsilon$, we get
  \begin{equation*}
  0\leq\frac{f\left(\phi(\rho)\right) F(\phi(\rho))\left(1-\frac 1\rho\right)}  { (1-\frac 1B+\sigma)(\rho-1)^{-\frac 1B+\sigma}} \leq \frac{{C}\left(1-\frac 1\rho\right)}{(1-\frac 1B+\sigma)\rho^{\frac 1B-\varepsilon} (\rho -1)^{-\frac 1B+\sigma}} \to 0 \ \ \ {\rm as}\ \ \ \rho\to \infty. 
   \end{equation*}
 
In order to establish the inequality \eqref{Stima2}, we first remark that $u'_{\infty}(r)\sim {\tilde{u}} '(r)$, as $r\to 0$, because
$$
\lim_{r\to 0}\frac{u'_{\infty}(r)}{{\tilde{u}} '(r)}=\lim_{\rho\to \infty}\frac{{\phi}'(\rho)+{\eta}'{(\rho)}}{{\phi}'(\rho)}=1
$$
thanks to the property $\lim_{\rho\to \infty}\frac{{\eta}'{(\rho)}}{{\phi}'(\rho)}=0$  (see  Lemma \ref{lemma:5.1}).
Therefore, the inequality \eqref{Stima2} is a consequence of
$$
\lim_{r\to 0}|{\tilde{u}}'(r)| r\left(-2\log r\right)^{\frac 1B-\sigma}=2\lim_{\rho\to \infty}|\phi'(\rho)|(\rho-1)^{\frac 1B-\sigma}=0,
$$ 
which can be proven 
in a similar way to the proof of \eqref{Stima1}.
\end{proof}

\begin{proof}[Proof of Proposition~\ref{proposition:2.3}]
Using the estimates  \eqref{Stima1}  and    \eqref{Stima2},  we can prove  that the function $u_\infty$ satisfies the differential equation in the sense of distributions in $B_{R}$.
 We use similar arguments as in \cite[Page 265]{BL}, page 265
and in \cite[Pages 261--262]{NS}. Let $\varphi$ be a
${ C^{\infty}}$ function with compact support in
$B_{R}$. We prove that
$$
\int_{B_R} {u_\infty}\ \Delta \varphi+\ f({u_\infty})\ \varphi\
dx=0.
$$
Indeed, let $\Phi(r)$ be a $ C^{\infty}(\mathbb{R})$
function, $0\leq \Phi(r)\leq 1$ such that
\begin{equation*}
\Phi (r)=\left\{\begin{split} 1\ \ \ &\ \text{if}\ \ r<1/2,\\
0\ \ \ &\ \text{if}\ \ r\geq 1,\\
\end{split}\right.
\end{equation*}
and
$\Phi_\varepsilon(|x|)=\Phi\left(\frac{\log|x|}{\log\varepsilon}\right)$
for any $x\neq 0$ (these cut-off functions are the
same as those used in \cite{BL}).  By a direct computation  for
small $\varepsilon>0$ we get $\Phi_\varepsilon (|x|)=1$ for
$|x|>\sqrt \varepsilon$ and $\Phi_\varepsilon (|x|)=0$ for
$|x|\leq \varepsilon$ and  for  $x\neq 0$, we get
$\Phi_\varepsilon(|x|)\to 1$ for  $\varepsilon \to 0^+$. By the
Dominated Convergence Theorem, since $ u_\infty$ 
belongs to $L^1(B_{R})$, we have
\begin{equation*}
\begin{aligned}
&\int_{B_{R}} u_\infty \, \Delta  \varphi+ f( u_\infty)\, \varphi\
dx\\&=\lim_{\varepsilon \to 0^+}\int_{B_{R}}\Phi_\varepsilon
  u_\infty  \, \Delta \varphi+\Phi_\varepsilon f( u_\infty)\ \varphi\ dx
\\&=\lim_{\varepsilon \to 0^+}\int_{B_{R}}\Phi_\varepsilon\,
 \Delta u_\infty\, \varphi\
dx+2\int_{B_{R}}\nabla\Phi_\varepsilon\cdot \nabla u_\infty \
\varphi\ dx+\int_{B_{R}}\Delta\Phi_\varepsilon \, u_\infty\,
\varphi\ dx
+\int_{B_{R}}\Phi_\varepsilon\, f(u_\infty)\, \varphi\ dx.\\
\end{aligned}
\end{equation*}
Since $ u_\infty$ is a classical solution of the elliptic
equation in $B_R\setminus\{0\}$ we obtain
\begin{equation*}
\begin{aligned}
&\lim_{\varepsilon \to 0^+}\int_{B_R}\Phi_\varepsilon\,  u_\infty\,
\Delta \varphi+\Phi_\varepsilon f( u_\infty)\ \varphi\ dx
\\&=\lim_{\varepsilon \to 0^+}2\int_{B_{R}}\nabla\Phi_\varepsilon\cdot
\nabla  u_\infty\ \varphi \
dx+\int_{B_{R}}\Delta\Phi_\varepsilon\,  u_\infty\, \varphi\
dx.\\
\end{aligned}
\end{equation*}
Since $$ \Delta \Phi_\varepsilon=\Phi^{''}\Big(\frac{\log r}{\log
\varepsilon}\Big)\frac{1}{r^2(\log\varepsilon)^2}$$
thanks to \eqref{Stima1},
for $0<\sigma<\frac{1}{B}$, we have
 $$
 \Big|\int_{B_{R}}  u_\infty \, \Delta \Phi_\varepsilon \, \varphi\ dx\Big|\leq
 \frac{C}{(\log\varepsilon)^2}\int_{\varepsilon}^{\sqrt{\varepsilon}}\frac{\left(-2\log r\right)^{1-\frac 1B+\sigma}}{r}\ dr
 $$
 and
 $$
 \begin{aligned}
\displaystyle  \lim_{\varepsilon\to 0^+}
\frac{1}{(\log\varepsilon)^2}
\int_\varepsilon^{\sqrt{\varepsilon}}
\frac{\left(-2\log r\right)^{1-\frac 1B+\sigma}}{r}\ 
 dr
 =
 \lim_{\varepsilon \to 0}
 \frac{C}{
 \left(-\log \varepsilon\right)^{\frac 1B-\sigma}}=0,
\end{aligned}
$$
for some positive constant $C$.
In a similar way, using \eqref{Stima2}, we have
$$
 \Big|\int_{B_{R}} \nabla  u_\infty\cdot \nabla \Phi_\varepsilon\ \varphi \ dx \Big|\leq
 \frac{C}{(-\log\varepsilon)}\int_{\varepsilon}^{\sqrt{\varepsilon}}\frac1{r\left(-2\log r \right)^{\frac 1B-\sigma}}\ dr
 $$
and
$$
 \lim_{\varepsilon \to 0^+} \frac{\int_{\varepsilon}^{\sqrt{\varepsilon}}\frac{ dr}{r\left(-2\log r \right)^{\frac 1B-\sigma}
 }}{{(-\log\varepsilon)}}
 =
\lim_{\varepsilon \to 0^+} \frac{C}{\left(-\log
\varepsilon\right)^{\frac1B-\sigma}}=0.
$$
This proves that the function $u_\infty$ satisfies 
$-\Delta u_{\infty}=f(u_{\infty})$ in $B_R$
in the sense of distributions.
\end{proof}

\section{Examples of singular solutions}
\label{section:4}
In this section, we apply Theorem~\ref{theorem:1.2} to several examples.
The following lemma is useful to compute concrete examples
since we shall consider nonlinear terms in the form of $f(u)=e^{a(u)}$.

\begin{lemma}
\label{lemma:3.1}
Let $a\in C^{5}\bigl((0,\infty)\bigr)$ be a function satisfying 
$
\lim_{s\to \infty}a(s)
=
\lim_{s\to \infty}a'(s)
=\infty
$
and define $a_n\ (n=1,2,3,4,5)$ by
\[
a_1(s)=\frac{1}{a'(s)},\qquad a_{n+1}(s)=\frac{a_{n}'(s)}{a'(s)}\ \ (n=1,2,3,4).
\]
Suppose that 
\begin{equation}
\label{eq:3.1a}
\lim_{s\to \infty}a_n(s)e^{-a(s)}=0
\ \ 
\text{and}\ \ 
\int_s^{\infty}
|a_n'(t)|e^{-a(t)}dt<\infty
\ \ \text{for\ all\ }n=1,2,3,4.
\end{equation}
Furthermore, for $n=1,2,3,4$
we either assume 
$a_n= 0$ on $[s,\infty)$ for some $s>0$,
or
$a_n\neq 0$ on $[s,\infty)$ for some $s>0$
and
\begin{equation}
\label{eq:3.0}
\lim_{s\to \infty}\frac{a_{n+1}(s)}{a_n(s)}=0.
\end{equation}%
Then, for the function $f(s)=e^{a(s)}$ there hold
\begin{align}
&
\label{eq:3.2}
F(s)
=
\Bigl(a_1(s)+a_2(s)+a_3(s)+a_4(s)\bigl(1+o(1)\bigr)\Bigr)e^{-a(s)},
\\
&
\label{eq:3.5c}
\log F(s)
=
-a(s)-\log a'(s)+O\left(\frac{a''(s)}{a'(s)^2}\right),
\\
&
\label{eq:3.6b}
f(s)F(s)
=
\frac{1}{a'(s)}
\Bigl(1+o(1)\Bigr),
\\
&
\label{eq:3.7b}
f'(s)F(s)
=1-\frac{a''(s)}{a'(s)^2}+\frac{3a''(s)^2-a'(s)a'''(s)}{a'(s)^4}\Bigl(1+o(1)\Bigr),
\\
&
\label{eq:3.8b}
\frac{f(s)f''(s)F(s)}{f'(s)}
=
1
+\frac{2a''(s)^2-a'(s)a'''(s)}{a'(s)^4}
+
O\left(
\left|\frac{a''a'''}{a'^5}\right|
+
\left|\frac{a''''}{a'^4}\right|
+
\left|\frac{a''^3}{a'^6}\right|
\right),
\end{align}
as $s\to \infty$.
\end{lemma}
\begin{proof}[Proof of Lemma~\ref{lemma:3.1}]
We write $a=a(s)$ for the sake of simplicity.
Clearly it holds
\begin{equation}\label{eq:3.1}
f(s)=e^{a},\quad f'(s)=a'e^{a},\quad f''(s)=\Bigl(a''+a'^2\Bigr)e^{a}.
\end{equation}
We first give the proof assuming that
$a_n\neq 0$ on $[s,\infty)$ for some $s>0$
with \eqref{eq:3.0} for all $n=1,2,3,4$.
It follows from integration by parts with the aid of the assumption~\eqref{eq:3.1a} that
\[
F(s)
=
\int_s^{\infty} \frac{dt}{e^a}
=
\int_s^{\infty}\frac{1}{-a'}[e^{-a}]'{dt}
=
a_1e^{-a}+\int_s^{\infty}a_1' e^{-a}{dt}.
\]
Repeating this argument, we obtain from the definition of $a_n$ that
\[
F(s)
=
\Bigl(a_1+a_2+a_3+a_4\Bigr)e^{-a}
+
\int_s^{\infty}
a_4'e^{-a}{dt}.
\]
The assumption \eqref{eq:3.0} and de l'Hospital rule yield that
\[
\lim_{s\to \infty}
\frac{
\int_s^{\infty}
a_4'
e^{-a}{dt}
}
{
a_4 e^{-a}
}
=
\lim_{s\to \infty}
\frac{-a_4' e^{-a}}{a_4' e^{-a}-a'a_4 e^{-a}}
=
\lim_{s\to \infty}
\frac{\frac{a_4'}{a'a_4}}{1-\frac{a_4'}{a'a_4}}
=
\lim_{s\to \infty}
\frac{\frac{a_5}{a_4}}{1-\frac{a_5}{a_4}}
=
0.
\]
This 
yields
\eqref{eq:3.2}.
Remark that it holds
\begin{equation}\label{eq:3.0aaa}
F(s)
=
a_1\biggl(1+O\Big(\frac{a_2}{a_1}\Big)\biggr)e^{-a(s)}
=
a_1\Bigl(1+\frac{a_2}{a_1}+\frac{a_3}{a_1}(1+o(1))\Bigr)e^{-a(s)}
\end{equation}
by \eqref{eq:3.0}.
Then we obtain 
\eqref{eq:3.5c}, \eqref{eq:3.6b}, and \eqref{eq:3.7b}
from \eqref{eq:3.1}, \eqref{eq:3.0aaa}, 
and the fact
\[
a_1=\frac{1}{a'},\ \ 
a_2=-\frac{a''}{a'^3},\ \ 
a_3=\frac{-a'a'''+3a''^2}{a'^5}.
\]
It remains to compute 
$\frac{f(s)f''(s)F(s)}{f'(s)}$.
Applying \eqref{eq:3.1} and \eqref{eq:3.2},
we have
\[
\begin{aligned}
\frac{ff''F}{f'}
&
=
\frac{a''+a'^2}{a'}
\left(
a_1+a_2+a_3+a_4(1+o(1))
\right)
=
\frac{a''+a'^2}{a'^2}
\left(
1+a_1'+a_2'+a_3'(1+o(1))
\right)
\\
&
=
\Bigl(1-a_1'\Bigr)
\left(
1+a_1'+a_2'+a_3'(1+o(1))
\right)
=
1+(a_2'-a_1'^2)+(a_3'-a_1'a_2')+a_3'o(1).
\end{aligned}
\]
In the final equality we applied $O(a_1'a_3')=o(a_3')$ which follows from
$a_1'=a'a_2=\frac{a_2}{a_1}=o(1)$ by \eqref{eq:3.0}.
The equation \eqref{eq:3.8b} follows from
\[
\begin{aligned}
&a_2'-a_1'^2=\frac{2a''(s)^2-a'(s)a'''(s)}{a'(s)^4},
\ \ 
a_3'=\frac{-a'(s)^2a''''(s)+10a'(s)a''(s)a'''(s)-15a''(s)^3}{a'(s)^6},
\\
&
a_3'-a_1'a_2'
=
\frac{-a'(s)^2a''''(s)+9a'(s)a''(s)a'''(s)-12a''(s)^3}{a'(s)^6}.
\end{aligned}
\]

If $n\in \{1,2,3,4\}$ satisfies $a_n=0$ on $[s,\infty)$ for some $s>0$,
then the expansion of $F(s)$ in \eqref{eq:3.2} holds 
with $a_k= 0$ for $n\le k \le 4$.
Thus the same conclusion holds.
\end{proof}

\begin{example}\label{example:3.1}
Let $q>1$, $r\in \mathbb{R}$, and $f(u)=u^re^{u^q}$.
Then there exists $R>0$ and a singular solution 
${u}\in C^{\infty}(B_R\setminus\{0\})$ to \eqref{eq:1.1}
satisfying
\[
{u(x)}
=
\bigg(\log\frac{1}{|x|^2}-\frac{2q+r-1}{q}\log\log\frac{1}{|x|^2}+\log\frac{4(q-1)}{q^2}
\bigg)^{\frac{1}{q}}
+O\bigg(\frac{\log(-\log|x|)}{(-\log|x|)^{2-\frac{1}{q}}}\bigg)
\]
as $x\to 0$.
\end{example}
\begin{proof}
We start from computing the order of the singular function $\tilde u:=F^{-1}[G(v)]$.
To this end, we prove that
\begin{equation}\label{eq:3.4}
F^{-1}(w)
=
\left(
-\log w
-\frac{q+r-1}{q}\log(-\log w)
+{\log\frac{1}{q}}
\right)^{\frac{1}{q}}
+
O\bigg(\frac{\log(-\log w)}{(-\log w)^{2-\frac{1}{q}}}\bigg).
\end{equation}
Let $w=F(s)$. 
Since $f(s)=s^re^{s^q}=e^{s^q+r\log s}$, we apply Lemma~\ref{lemma:3.1} to $a(s)=s^q+r\log s$.
The function $a(s)$ satisfies all assumptions in Lemma~\ref{lemma:3.1}.
Then by 
\eqref{eq:3.5c}
it holds
\[
\log w=-s^q-(q+r-1)\log s+\log \frac{1}{q}+O(s^{-q}),
\]
and hence
\[
s=\left(-\log w-(q+r-1)\log s+\log \frac{1}{q}+O(s^{-q})\right)^{\frac{1}{q}}.
\]
It follows from $e^{-2s^q}\le w \le e^{-s^q/2}$ that 
\[
2s^q \ge -\log w \ge {\frac{s^q}{2}}.
\]
This yields
\[
\begin{aligned}
s
=
\left(
-\log w-(q+r-1)\log s
+\log \frac{1}{q}+O\left(\frac{1}{-\log w}\right)
\right)^{\frac{1}{q}}.
\end{aligned}
\]
Substituting $s$ into $\log s$, we have
\[
\begin{aligned}
\log s
&
= \frac{1}{q}\log 
\left[
-\log w-(q+r-1)\log s+\log \frac{1}{q}+O\left(\frac{1}{-\log w}\right)
\right]
\\
&
=
\frac{1}{q}\log(-\log w)+O\left(\frac{\log(-\log w)}{\log w}\right).
\end{aligned}
\]
All together, we have
\[
\begin{aligned}
s
&
=
\left(
-\log w
-\frac{q+r-1}{q}\log(-\log w)
+\log \frac{1}{q}
+O\left(\frac{\log(-\log w)}{\log w}
\right)
\right)^{\frac{1}{q}}
\\
&
=
\left(
-\log w
-\frac{q+r-1}{q}\log(-\log w)
+\log \frac{1}{q}
\right)^{\frac{1}{q}}
+
O\bigg(\frac{\log(-\log w)}{(-\log w)^{2-\frac{1}{q}}}\bigg).
\end{aligned}
\]
This proves \eqref{eq:3.4}.
Thanks to \eqref{eq:3.4}, we obtain
\begin{equation}
\label{eq:3.6}
\begin{aligned}
\tilde u(x)
&
=
F^{-1}[G(v)]
=
F^{-1}
\left[
\frac{q'}{4}
|x|^2 \left(\log \frac{1}{|x|^2}+1\right)
\right]
\\
&
=
\left(
\log\frac{1}{|x|^2}-\frac{2q+r-1}{q}\log\log\frac{1}{|x|^2}+\log\frac{4(q-1)}{q^2}
\right)^{\frac{1}{q}}
+O\bigg(\frac{\log(-\log|x|)}{(-\log|x|)^{2-\frac{1}{q}}}\bigg).
\end{aligned}
\end{equation}

We next show that the nonlinear term $f(s)=s^re^{s^{q}}$ satisfies all assumptions in Theorem~\ref{theorem:1.2}.
It is easy to check from Lemma~\ref{lemma:3.1} with $a(s)=s^q+r\log s$ that
\[
\begin{aligned}
&
\log F(s)
=
-a(s)-\log a'(s)+O\left(\frac{a''(s)}{a'(s)^2}\right)
=
-s^q-(q+r-1)\log s 
+\log\frac{1}{q}
+O\left(\frac{1}{s^{q}}\right),
\\
&
f'(s)F(s)
=1-\frac{a''(s)}{a'(s)^2}+\frac{3a''(s)^2-a'(s)a'''(s)}{a'(s)^4}\Bigl(1+o(1)\Bigr)
=
1
-
\left(1-\frac{1}{q}\right)\frac{1}{s^q}+O\left(\frac{1}{s^{2q}}\right),
\\
&\frac{f(s)f''(s)F(s)}{f'(s)}
\\
&
=
1+\frac{2a''(s)^2-a'(s)a'''(s)}{a'(s)^4}
+O\left(
\left|\frac{a''a'''}{a'^5}\right|
+
\left|\frac{a''''}{a'^4}\right|
+
\left|\frac{a''^3}{a'^6}\right|
\right)
\\
&=
1+\frac{2\Bigl(q(q-1)s^{q-2}-rs^{-2}\Bigr)^2-(qs^{q-1}+rs^{-1})\Bigl(q(q-1)(q-2)s^{q-3}+2rs^{-3}\Bigr)}{(qs^{q-1}+rs^{-1})^4}
+O\left(\frac{1}{s^{3q}}\right)
\\
&
=
1
+
\left(1-\frac{1}{q}\right)\frac{1}{s^{2q}}
+
O\left(\frac{1}{s^{3q}}\right).
\end{aligned}
\]
These equalities yield that
\begin{equation}\label{eq:3.5a}
\frac{1}{B_1[f](s)}=\log F(s)(f'(s)F(s)-1)=1-\frac{1}{q}+O\left(\frac{\log s}{s^q}\right)
\end{equation}
and
\begin{equation}\label{eq:3.5b}
\begin{aligned}
\frac{1}{B_2[f](s)}
&
=(\log F)^2f'F\left(\frac{ff''F}{f'}-1\right)
\\
&
=
\Bigl(s^{2q}+2(q+r-1)s^q\log s 
+O\left({s^{q}}\right)\Bigr)
\\
&\qquad \times
\left(1-\left(1-\frac{1}{q}\right)\frac{1}{s^q}+O\left(\frac{1}{s^{2q}}\right)\right)
\times 
\left(\left(1-\frac{1}{q}\right)\frac{1}{s^{2q}}
+
O\left(\frac{1}{s^{3q}}\right)\right)
\\
&
=
\left(1-\frac{1}{q}\right)+O\left(\frac{\log s}{s^q}\right).
\end{aligned}
\end{equation}
Clearly 
$g(s)=\frac{1}{qq'}\frac{e^{s^q}}{s^{2q-1}}$
also satisfies 
\begin{equation}\label{eq:1.20}
\frac{1}{B_1[g](s)}=1-\frac{1}{q}+O\left(\frac{\log s}{s^{q}}\right),
\qquad 
\frac{1}{B_2[g](s)}=1-\frac{1}{q}+O\left(\frac{\log s}{s^{q}}\right).
\end{equation}
Combining \eqref{eq:3.6}, \eqref{eq:3.5a}, \eqref{eq:3.5b}, \eqref{eq:1.20}, and \eqref{eq:1.3a} {with $B=q'$},  we have
\[
R_1=R_2=
O\left(\frac{\log \tilde u}{\tilde u^q}\right)
+
O\left(\frac{\log v}{v^q}\right)
=
O\left(\frac{\log(-\log |x|)}{-\log |x|}\right).
\]
This proves that $f(s)=s^re^{s^q}$ satisfies the assumption \eqref{eq:1.1a}.
Furthermore,
since 
\[
f(s)F(s)=O\left(\frac{1}{(s^q+r\log s)'}\right)=O\left(s^{-(q-1)}\right), 
\]
it holds
\[
f\bigl(\tilde u(x)\bigr)F\bigl(\tilde u(x)\bigr)\sup_{|y|\le |x|}
\bigl(R_1(y)+R_2(y)\bigr)
=
O\left(\frac{\log (-\log |x|)}{(-\log |x|)^{2-\frac{1}{q}}}\right).
\]
This together with \eqref{eq:3.6} yields that the singular solution ${u}$ constructed in 
Theorem~\ref{theorem:1.2}
satisfies
\[
\begin{aligned}
{u(x)}
&
=
\tilde u(x)+O\left(\frac{\log (-\log |x|)}{(-\log |x|)^{2-\frac{1}{q}}}\right)
\\
&
=
\left(
\log\frac{1}{|x|^2}-\frac{2q+r-1}{q}\log\log\frac{1}{|x|^2}+\log\frac{4(q-1)}{q^2}
\right)^{\frac{1}{q}}
+
O\left(\frac{\log(-\log|x|)}{(-\log|x|)^{2-\frac{1}{q}}}\right).
\end{aligned}
\]
This completes the proof of Example~\ref{example:3.1}.
\end{proof}

\begin{example}\label{example:3.2}
Let $q>1$, 
$\frac{q}{2}>r  >0$,
and $f(u)=e^{u^q+u^r}$.
Then there exists $R>0$ and a singular solution 
${u}\in C^{\infty}(B_R\setminus\{0\})$ to \eqref{eq:1.1}
satisfying
\[
\begin{aligned}
{u(x)}
=
\Biggl[
-\log w
&
-
\biggl\{
-\log w
+\sum_{k=0}^{n-2}\left(\frac{r}{q}\right)^k{(-\log w)^{(k+1)\frac{r}{q}-k}}
-\frac{q-1}{q}\log(-\log w)
\biggr\}^{\frac{r}{q}}
\\
&
-\frac{q-1}{q}\log(-\log w) +\log \frac{1}{q}
\Biggr]^{\frac{1}{q}}
+O\left(\frac{1}{(-\log w)^{2-\frac{1}{q}-\frac{r}{q}}}\right)
\end{aligned}\]
as $x\to 0$,
where {$n\in \mathbb{N}$ with $n\ge 2$ is  such that
\[
\frac{1}{1-\frac{r}{q}}<n\le \frac{1}{1-\frac{r}{q}}+1,
\]
and }
\[
w
:=\frac{q'}{4}|x|^2\left(\log \frac{1}{|x|^2}+1\right).
\]
\end{example}

\begin{remark}
Higher expansion in ${u}$ is needed to obtain the expected decay order 
$O\left(\frac{1}{(-\log w)^{2-\frac{1}{q}-\frac{r}{q}}}\right)$.
Indeed, a more simple expression is possible:
\[
{u(x)}
=
\left[
-\log w
-
(-\log w)^{\frac{r}{q}}-\frac{q-1}{q}\log(-\log w)+\log \frac{1}{q}
\right]^{\frac{1}{q}}
+O\left(\frac{1}{(-\log w)^{3-\frac{1}{q}-3\frac{q}{r}}}\right).
\]
Then this decay order is worse than the order appearing from the remainder term of the singular solution
$fF(R_1+R_2)=O\left(\frac{1}{(-\log w)^{2-\frac{1}{q}-\frac{r}{q}}}\right)$.
By considering the higher expansion for  $F^{-1}(w)$, we get the same decay order from $fF(R_1+R_2)$.
\end{remark}

\begin{proof}
We first prove  that
\begin{equation}\label{eq:3.8}
\begin{aligned}
F^{-1}(w)
=
\Biggl[
-\log w
&
-
\biggl\{
-\log w
+\sum_{k=0}^{n-2}\left(\frac{r}{q}\right)^k{(-\log w)^{(k+1)\frac{r}{q}-k}}
-\frac{q-1}{q}\log(-\log w)
\biggr\}^{\frac{r}{q}}
\\
&
-\frac{q-1}{q}\log(-\log w) +\log \frac{1}{q}
\Biggr]^{\frac{1}{q}}
+O\left(\frac{1}{(-\log w)^{2-\frac{1}{q}-\frac{r}{q}}}\right),
\end{aligned}
\end{equation}
by the similar argument as in Example~\ref{example:3.1}.
Let $w=F(s)$.
We apply Lemma~\ref{lemma:3.1} to $a(s)=s^q+s^r$.
The function $a(s)$ satisfies all assumptions in Lemma~\ref{lemma:3.1}.
Then, by \eqref{eq:3.2} we have
\[
F(s)=\frac{e^{-(s^q+s^r)}}{qs^{q-1}+rs^{r-1}}+O\left(\frac{e^{-(s^q+s^r)}}{s^{2q-1}}\right).
\]
Hence it holds
\begin{equation}\label{eq:3.7}
\begin{aligned}
\log w
&
=
-s^q-s^r-\log (qs^{q-1}+rs^{r-1})+O\left(\frac{1}{s^q}\right)
\\
&
=
-s^q-s^r-(q-1)\log s +\log \frac{1}{q}+O\left(\frac{1}{s^{q-r}}\right).
\end{aligned}
\end{equation}
This shows us that
\begin{equation}\label{eq:3.8a}
\begin{aligned}
s^r
&
=
(s^{q})^{\frac{r}{q}}
=
\left(
-\log w
-s^r-(q-1)\log s +\log \frac{1}{q}+O\left(\frac{1}{s^{q-r}}\right)
\right)^{\frac{r}{q}}
\\
&
=
\left(
-\log w
-s^r-(q-1)\log s +\log \frac{1}{q}+O\left(\frac{1}{(-\log w)^{1-\frac{r}{q}}}\right)
\right)^{\frac{r}{q}},
\end{aligned}
\end{equation}
where 
we used
$e^{-2s^q}\le w\le e^{-s^q}$
and hence 
$ O(s^{-(q-r)})=O((-\log w)^{-(1-\frac{r}{q})})$.
Moreover, 
since 
$s^r=O\Bigl((-\log w)^{\frac{r}{q}}\Bigr)$, we have
\[
s^r
=
\left(-\log w+O\Bigl((-\log w)^{\frac{r}{q}}\Bigr)\right)^{\frac{r}{q}}
=
(-\log w)^{\frac{r}{q}}+O\left(\frac{1}{(-\log w)^{1-\frac{2r}{q}}}\right).
\]
This yields
\[
r\log s
=
\log
\left[
(-\log w)^{\frac{r}{q}}+O\left(\frac{1}{(-\log w)^{1-\frac{2r}{q}}}\right)
\right]
=
\frac{r}{q}\log (-\log w)
+O\left(\frac{1}{(-\log w)^{1-\frac{r}{q}}}\right).
\]
Substituting $\log s$ into \eqref{eq:3.8a}, we obtain
\begin{equation}\label{eq:3.9a}
s^r
=
\left[
-\log w  -s^r -\frac{q-1}{q}\log (-\log w)+\log \frac{1}{q}
+O\left(\frac{1}{(-\log w)^{1-\frac{r}{q}}}\right)
\right]^{\frac{r}{q}}.
\end{equation}
We determine the order of $s^r$ from the above equality~\eqref{eq:3.9a}.
Let
\[
a=-\log w,\quad b=s^r,\quad c=-\frac{q-1}{q}\log(-\log w)+\log \frac{1}{q},\quad d=O\left(\frac{1}{(-\log w)^{1-\frac{r}{q}}}\right)
\]
for the sake of simplicity.
Then \eqref{eq:3.9a} becomes
\[
b=(a+b+c+d)^{\frac{r}{q}}.
\]
Substituting $b$ into the right hand side, we have
\[
\begin{aligned}
b
&
=
\Biggl[
a+(a+b+c+d)^{\frac{r}{q}}+c+d
\Biggr]^{\frac{r}{q}}
\\
&
=
\Biggl[
a+a^{\frac{r}{q}}+\frac{r}{q}\frac{b+c+d}{a^{1-\frac{r}{q}}}+O\left(\frac{{(b+c+d)^2}}{a^{2-\frac{r}{q}}}\right)+c+d
\Biggr]^{\frac{r}{q}}
\\
&
=
\Biggl[
a+a^{\frac{r}{q}}+\frac{r}{q}\frac{1}{a^{1-\frac{r}{q}}}b+\frac{r}{q}\frac{c}{a^{1-\frac{r}{q}}}+c+d
\Biggr]^{\frac{r}{q}}.
\end{aligned}
\]
Again, substituting $b$ into the right hand side, we have
\[
\begin{aligned}
b
&
=
\Biggl[
a+a^{\frac{r}{q}}+\frac{r}{q}\frac{1}{a^{1-\frac{r}{q}}}
\Biggl(
a+a^{\frac{r}{q}}+\frac{r}{q}\frac{1}{a^{1-\frac{r}{q}}}b+\frac{r}{q}\frac{c}{a^{1-\frac{r}{q}}}+c+d
\Biggr)^{\frac{r}{q}}
+\frac{r}{q}\frac{c}{a^{1-\frac{r}{q}}}+c+d
\Biggr]^{\frac{r}{q}}
\\
&
=
\Biggl[
a+a^{\frac{r}{q}}+\frac{r}{q}\frac{1}{a^{1-2\frac{r}{q}}}
+\left(\frac{r}{q}\right)^2\frac{1}{a^{2-3\frac{r}{q}}}
+\left(\frac{r}{q}\right)^3\frac{b}{a^{3-3\frac{r}{q}}}
+\frac{r}{q}\frac{c}{a^{1-\frac{r}{q}}}+c+d
\Biggr]^{\frac{r}{q}}.
\end{aligned}
\]
{Since $2r<q$ there exists} $n\in \mathbb{N}$ with $n\ge 2$ such that
\[
\frac{1}{1-\frac{r}{q}}<n\le \frac{1}{1-\frac{r}{q}}+1,
\]
namely $(n-1)-n\frac{r}{q}\le 1-\frac{r}{q}<n-(n+1)\frac{r}{q}$.
For this $n$, it holds $k\frac{r}{q}-(k-1)>0$ for $k=1,\ldots, n-1$ and $n\frac{r}{q}-(n-1)\le 0$.
Repeating the above argument $n$ times, we obtain
\[
b
=
\Biggl[
a
+
\sum_{k=0}^{n-2}
\left(\frac{r}{q}\right)^ka^{(k+1)\frac{r}{q}-k}
+
\left(\frac{r}{q}\right)^{n-1}
\frac{1}
{a^{n-1-n\frac{r}{q}}}
+
\left(\frac{r}{q}\right)^n
{b}{a^{n\frac{r}{q}-n}}
+\frac{r}{q}\frac{c}{a^{1-\frac{r}{q}}}+c+d
\Biggr]^{\frac{r}{q}}.
\]
It follows from the definitions of $a,b,c,d$ and
$n$ 
that 
\[
\left(\frac{r}{q}\right)^{n-1}
\frac{1}
{a^{n-1-n\frac{r}{q}}}
+
\left(\frac{r}{q}\right)^n
{b}{a^{n\frac{r}{q}-n}}
+\frac{r}{q}\frac{c}{a^{1-\frac{r}{q}}}+d
=
O\left(
\frac{1}
{a^{n-1-n\frac{r}{q}}}
\right)
=
O\left(
\frac{1}{(-\log w)^{n-1-n\frac{r}{q}}}
\right).
\]
Therefore,
\begin{equation}\label{eq:3.10a}
\begin{aligned}
s^r
&
=
\Biggl\{
-\log w
+
\sum_{k=0}^{n-2}\left(\frac{r}{q}\right)^k{(-\log w)^{(k+1)\frac{r}{q}-k}}
\\
&
\qquad \qquad 
-
\frac{q-1}{q}\log(-\log w)
+
\log \frac{1}{q}
+
O\left(
\frac{1}{(-\log w)^{n-1-n\frac{r}{q}}}
\right)
\Biggr\}^{\frac{r}{q}}
\\
&
=
\Biggl\{
-\log w
+\sum_{k=0}^{n-2}\left(\frac{r}{q}\right)^k{(-\log w)^{(k+1)\frac{r}{q}-k}}
-
\frac{q-1}{q}\log(-\log w)
\Biggr\}^{\frac{r}{q}}
+
O\left(
\frac{1}{(-\log w)^{n-(n+1)\frac{r}{q}}}
\right).
\end{aligned}
\end{equation}
In conclusion, by \eqref{eq:3.7} and \eqref{eq:3.10a}
it holds
\[
\begin{aligned}
s^q
&
=
-\log w
-
\Biggl\{
-\log w
+\sum_{k=0}^{n-2}\left(\frac{r}{q}\right)^k{(-\log w)^{(k+1)\frac{r}{q}-k}}
-
\frac{q-1}{q}\log(-\log w)
\Biggr\}^{\frac{r}{q}}
\\
&
\qquad
-\frac{q-1}{q}\log(-\log w) +\log \frac{1}{q}+O\left(\frac{1}{(-\log w)^{1-\frac{r}{q}}}\right),
\end{aligned}
\]
where we used 
$n-(n+1)\frac{r}{q}\ge 1-\frac{r}{q}$.
This yields
\[
\begin{aligned}
s
=
\Biggl[
-\log w
&
-
\biggl\{
-\log w
+\sum_{k=0}^{n-2}\left(\frac{r}{q}\right)^k{(-\log w)^{(k+1)\frac{r}{q}-k}}
-\frac{q-1}{q}\log(-\log w)
\biggr\}^{\frac{r}{q}}
\\
&
-\frac{q-1}{q}\log(-\log w) +\log \frac{1}{q}
\Biggr]^{\frac{1}{q}}
+O\bigg(\frac{1}{(-\log w)^{2-\frac{1}{q}-\frac{r}{q}}}\bigg),
\end{aligned}
\]
and hence \eqref{eq:3.8} is proved.
As a consequence, we obtain 
the singularity of the function 
$\tilde u:=F^{-1}(w)$
by taking
$
w
:=
G(v)=\frac{q'}{4}|x|^2\left(\log \frac{1}{|x|^2}+1\right)
$.

We next check that $f(s)=e^{s^q+s^r}$ with $q>1$ and 
$\frac{q}{2}>r>0$
 satisfies all the assumptions in Theorem~\ref{theorem:1.2}.
Applying Lemma~\ref{lemma:3.1} with $a(s)=s^q+s^r$, we have 
\[
\begin{aligned}
\log F(s)
&=
-a(s)-\log a'(s)+O\left(\frac{a''(s)}{a'(s)^2}\right)
=
-s^q-s^r
+O\left(\log s\right),
\\
f'(s)F(s)
&
=
1-\frac{a''(s)}{a'(s)^2}+\frac{3a''(s)^2-a'(s)a'''(s)}{a'(s)^4}\Bigl(1+o(1)\Bigr)
\\
&
=
1-\frac{q(q-1)s^{q-2}+r(r-1)s^{r-2}}{(qs^{q-1}+rs^{r-1})^2}+O\left(\frac{1}{s^{2q}}\right)
\\
&
=
1
-
\left(1-\frac{1}{q}\right)\frac{1}{s^q}+O\left(\frac{1}{s^{2q-r}}\right),
\\
\frac{f(s)f''(s)F(s)}{f'(s)}
&
=
1+
\frac{2a''(s)^2-a'(s)a'''(s)}{a'(s)^4}
+
O\left(
\left|\frac{a''a'''}{a'^5}\right|
+
\left|\frac{a''''}{a'^4}\right|
+
\left|\frac{a''^3}{a'^6}\right|
\right).
\end{aligned}
\]
One can compute easily that
\[
\begin{aligned}
&
\frac{2a''(s)^2-a'(s)a'''(s)}{a'(s)^4}
\\
&
=
\frac{
2\Bigl(q(q-1)s^{q-2}+r(r-1)s^{r-2}\Bigr)^2-
(qs^{q-1}+rs^{r-1})
\Bigl(q(q-1)(q-2)s^{q-3}+r(r-1)(r-2)s^{r-3}\Bigr)
}
{
(qs^{q-1}+rs^{r-1})^4
}
\\
&
=
\left(1-\frac{1}{q}\right)\frac{1}{s^{2q}}
+
O\left(\frac{1}{s^{3q-r}}\right)
\end{aligned}
\]
and
\[
O\left(
\left|\frac{a''a'''}{a'^5}\right|
+
\left|\frac{a''''}{a'^4}\right|
+
\left|\frac{a''^3}{a'^6}\right|
\right)
=
O\left(\frac{1}{s^{3q}}\right).
\]
Then it holds
\begin{equation}\label{eq:3.5cc}
\frac{1}{B_1[f](s)}=\log F(s)(f'(s)F(s)-1)=1-\frac{1}{q}+O\left(\frac{1}{s^{q-r}}\right)
\end{equation}
and
\begin{equation}\label{eq:3.5d}
\begin{aligned}
\frac{1}{B_2[f](s)}
&
=(\log F)^2f'F\left(\frac{ff''F}{f'}-1\right)
\\
&
=
\Bigl(s^{2q}+2s^{q+r}
+{O\left(s^{q}\log s\right)}
\Bigr)
\\
&\qquad \times
\left(1-\left(1-\frac{1}{q}\right)\frac{1}{s^q}+O\left(\frac{1}{s^{2q-r}}\right)\right)
\times 
\left(
\left(1-\frac{1}{q}\right)\frac{1}{s^{2q}}
+
O\left(\frac{1}{s^{3q-r}}\right)\right)
\\
&
=
\left(1-\frac{1}{q}\right)+O\left(\frac{1}{s^{q-r}}\right).
\end{aligned}
\end{equation}
Combining \eqref{eq:3.5cc}, \eqref{eq:3.5d}, and \eqref{eq:1.20}, we obtain
\[
\begin{aligned}
R_1=R_2
&
=
O\left(\frac{1}{\tilde u^{q-r}}\right)
+
O\left(\frac{\log v}{v^{q}}\right)
\\
&
=
O\left(\frac{1}{(-\log |x|)^{1-\frac{r}{q}}}\right)
+
O\left(\frac{\log(-\log |x|)}{-\log x}\right)
=
O\left(\frac{1}{(-\log |x|)^{1-\frac{r}{q}}}\right).
\end{aligned}
\]
Since $1-\frac{r}{q}>\frac{1}{2}$ by $\frac{q}{2}>r$, 
the assumption \eqref{eq:1.1a} is satisfied.
Furthermore, it follows from
\[
f(s)F(s)=O\left(\frac{1}{(s^q+s^r)'}\right)=O\left(s^{-(q-1)}\right)
\]
that
\begin{equation}\label{eq:3.11a}
f\bigl(\tilde u(x)\bigr)F\bigl(\tilde u(x)\bigr)\sup_{|y|\le |x|}
\bigl(R_1(y)+R_2(y)\bigr)
=
O\left(\frac{1}{(-\log |x|)^{2-\frac{1}{q}-\frac{r}{q}}}\right).
\end{equation}
Therefore, \eqref{eq:3.8} and \eqref{eq:3.11a}
yield 
the conclusion of Example~\ref{example:3.2}.
\end{proof}


\begin{example}\label{example:3.3}
Let $q>1$, $r\in \mathbb{R}$, and $f(u)=e^{u^q(\log u)^r}$. Then, $B=q'$ but the assumption~\eqref{eq:1.1a} 
is not fulfilled.
\end{example}
\begin{proof}
Let $q>1$, $r\in \mathbb{R}$, and $a(s)=s^q(\log s)^r$.
The function $a(s)$ satisfies all the assumptions in Lemma~\ref{lemma:3.1}.
Applying Lemma~\ref{lemma:3.1} to $a(s)$, we have
\[
\begin{aligned}
\log F(s)
&
=
-a(s)-\log a'(s)+O\left(\frac{a''(s)}{a'(s)^2}\right)
=
-s^q(\log s)^r-{(q-1)\log s }+O\left(\log \log s\right),
\\
f'(s)F(s)
&
=1-\frac{a''(s)}{a'(s)^2}+\frac{3a''(s)^2-a'(s)a'''(s)}{a'(s)^4}\Bigl(1+o(1)\Bigr)
\\
&
=
1
-
\left(1-\frac{1}{q}\right)\frac{1}{s^q(\log s)^r}+O\left(\frac{1}{s^{q}(\log s)^{r+1}}\right),
\\
\frac{f(s)f''(s)F(s)}{f'(s)}
&
=
1+\frac{2a''(s)^2-a'(s)a'''(s)}{a'(s)^4}
+
O\left(
\left|\frac{a''a'''}{a'^5}\right|
+
\left|\frac{a''''}{a'^4}\right|
+
\left|\frac{a''^3}{a'^6}\right|
\right).
\end{aligned}
\]
Here we have
\[
\begin{aligned}
&
2a''(s)^2-a'(s)a'''(s)
\\
&
=
2\Bigl(q(q-1)s^{q-2}(\log s)^r+O\bigl(s^{q-2}(\log s)^{r-1}\bigr)\Bigr)^2
\\
&
\quad -
(qs^{q-1}(\log s)^r+O\bigl(s^{q-1}(\log s)^{r-1}\bigr))
\Bigl(q(q-1)(q-2)s^{q-3}(\log s)^r+O\bigl(s^{q-3}(\log s)^{r-1}\bigr)\Bigr)
\\
&
=
q^3(q-1)s^{2q-4}(\log s)^{2r}+O\Bigl(s^{2q-4}(\log s)^{2r-1}\Bigr)
\end{aligned}
\]
and hence
\[
\begin{aligned}
\frac{
2a''(s)^2-a'(s)a'''(s)
}{a'(s)^4}
&
=
\frac{
q^3(q-1)s^{2q-4}(\log s)^{2r}+O\Bigl(s^{2q-4}(\log s)^{2r-1}\Bigr)}
{
(qs^{q-1}(\log s)^r+O\bigl(s^{q-1}(\log s)^{r-1}\bigr))^4
}
\\
&
=
\left(1-\frac{1}{q}\right)\frac{1}{s^{2q}(\log s)^{2r}}
+
O\left(\frac{1}{s^{2q}(\log s)^{2r+1}}\right).
\end{aligned}
\]
Furthermore, it holds
\[
O\left(
\left|\frac{a''a'''}{a'^5}\right|
+
\left|\frac{a''''}{a'^4}\right|
+
\left|\frac{a''^3}{a'^6}\right|
\right)
=
O\left(\frac{1}{s^{3q}(\log s)^{3r}}\right).
\]
Then, we obtain
\begin{equation}\label{eq:3.5e}
\frac{1}{B_1[f](s)}=\log F(s)(f'(s)F(s)-1)=1-\frac{1}{q}+O\left(\frac{1}{\log s}\right)
\end{equation}
and
\begin{equation}\label{eq:3.5f}
\begin{aligned}
\frac{1}{B_2[f](s)}
&
=(\log F)^2f'F\left(\frac{ff''F}{f'}-1\right)
\\
&
=
\Bigl(s^{2q}(\log s)^{2r}+{2(q-1)s^q(\log s)^{r+1}} +O\bigl((\log s)^2\bigr)\Bigr)
\\
&\qquad \times
\left(
1
-
\left(1-\frac{1}{q}\right)\frac{1}{s^q(\log s)^r}+O\left(\frac{1}{s^{q}(\log s)^{r+1}}\right)\right)
\\
&\qquad \qquad 
\times 
\left(
\left(1-\frac{1}{q}\right)\frac{1}{s^{2q}(\log s)^{2r}}
+
O\left(\frac{1}{s^{2q}(\log s)^{2r+1}}\right)
\right)
\\
&
=
1-\frac{1}{q}+O\left(\frac{1}{\log s}\right).
\end{aligned}
\end{equation}
This shows us that $\frac{1}{B}=\lim_{s\to \infty}\frac{1}{B_2[f](s)}=1-\frac{1}{q}$.

It remains to show that $f(u)=e^{u^q(\log u)^r}$ does not satisfy the assumption \eqref{eq:1.1a}.
Fix small $\varepsilon>0$. 
By \eqref{eq:3.2},
we see{
\[
e^{-s^{q+\varepsilon}}
\le 
F(s)=\frac{e^{-s^q(\log s)^r}}{qs^{q-1}(\log s)^r+rs^{q-1}(\log s)^{r-1}}\left(1+O\left(\frac{1}{s^q(\log 
s)^r}\right)\right)
\le
e^{-s^{q-\varepsilon}}
\]
for sufficiently large $s>0$.
Therefore, it holds
\[
(-\log w)^{\frac{1}{q+\varepsilon}} \le F^{-1}(w)\le (-\log w)^{\frac{1}{q-\varepsilon}}
\]
for sufficiently small $w>0$. }
This yields 
\[
\log \tilde{u}
=
\log\left(F^{-1}\left[\frac{B}{4}|x|^2\left(\log \frac{1}{|x|^2}+1\right)\right]\right)
=
O(\log \log |x|)
\]
as $|x|\to \infty$. 
This, together with \eqref{eq:3.5e}, \eqref{eq:3.5f}, and \eqref{eq:1.20}, {implies that}
\[
R_1(x)=R_2(x)=
O\left(\frac{1}{\log \tilde u}\right)
+
O\left(\frac{\log v}{v^{q}}\right)
=
O\left(\frac{1}{\log \log |x|}\right)
+
O\left(\frac{\log \log |x|}{\log |x|}\right)
=
O\left(\frac{1}{\log \log |x|}\right)
\]
as $|x|\to \infty$. This proves that 
$f(u)=e^{u^q(\log u)^r}$ does not satisfy \eqref{eq:1.1a}.
\end{proof}

\begin{example}\label{example:3.4}
Let $q\ge 1$ and $f(u)=e^{e^{u^q}}$. Then, $B=1$ for all $q\ge 1$.
Furthermore, the assumption~\eqref{eq:1.1a} 
is fulfilled for $q=1$ while \eqref{eq:1.1a} is not fulfilled for $q>1$.
If $q=1$, then there exists $R>0$ and a singular solution 
${u}\in C^{\infty}(B_R\setminus\{0\})$ to \eqref{eq:1.1}
satisfying
\[
\begin{aligned}
{u(x)}
=
\log
\Big[
-\log w -\log(-\log w)
\Big]
+O\Big(\frac{
\log(-\log w)
}{\left(-\log w\right)^2}\Big)
\end{aligned}
\]
as $x\to 0$,
where
\[
w=\frac{1}{4}|x|^2\Big(\log \frac{1}{|x|^2}+1\Big).
\]
\end{example}

\begin{proof}
Let $q\ge 1$ and $f(s)=e^{e^{s^q}}$, hence $a(s)=e^{s^q}$.
Then $a(s)$ satisfies all assumptions in Lemma~\ref{lemma:3.1}.
We first show that
\begin{equation}\label{eq:example4.40}
F^{-1}(w)
=
\left(
\log
\Big[
 -\log w -\log(-\log w)-\frac{q-1}{q}\log\log(-\log w)+\log \frac{1}{q}
\Big]
\right)^{\frac{1}{q}}
+
O\bigg(
\frac{\big(\log (-\log w)\big)}{(-\log w)^2}
\bigg)
\end{equation}
for $w>0$.
Let $w=F(s)$. 
Since $a'(s)=qs^{q-1}e^{s^q}$,
we have
\[
a_1=\frac{1}{qs^{q-1}e^{s^q}},
\qquad
a_2
=
\dfrac{
\frac{-q(q-1)s^{q-2}e^{s^q}-q^2s^{2q-2}e^{s^q}}{q^2s^{2q-2}e^{2s^q}}
}
{
qs^{q-1}e^{s^q}
}
=
-\frac{q-1}{q^2}\frac{1}{s^{2q-1}e^{2s^q}}
-\frac{1}{q}\frac{1}{s^{q-1}e^{2s^q}}.
\]
Then, by Lemma~\ref{lemma:3.1} we have
\[
w=F(s)=\frac{1}{qs^{q-1}e^{s^q}}\left[1+O(e^{-s^q})\right]e^{-e^{s^q}},
\]
and hence
\begin{equation}\label{eq:example4.4}
\log w
=
-e^{s^q}
-s^q-(q-1)\log s +\log \frac{1}{q}+O(e^{-s^q}).
\end{equation}
Therefore, 
\begin{equation}\label{eq:example4.42}
s^q=\log \Bigl[-\log w-s^q-(q-1)\log s +\log \frac{1}{q}+O(e^{-s^q})\Bigr].
\end{equation}
It follows from \eqref{eq:example4.4}
that
$e^{s^q}\le -\log w \le e^{2s^q}$
and hence
\[s^q \le \log (-\log w)\le 2s^q.\]
This together with \eqref{eq:example4.42}
yields
\[
s^q=\log(-\log w)+O\left(\frac{\log(-\log w)}{-\log w}\right).
\]
By \eqref{eq:example4.42}, we have
\[
s^q
=
\log
\left[
 -\log w -\log(-\log w)-\frac{q-1}{q}\log\log(-\log w)+\log \frac{1}{q}
 +O\left(\frac{\log(-\log w)}{-\log w}\right)
\right].
\]
Since $(1+t)^{\frac{1}{q}}=1+O(t)\ (t\to 0)$, we obtain \eqref{eq:example4.40}.

We next show that \eqref{eq:1.1a} is fulfilled if $q=1$ but it is not fulfilled if $q>1$.
Since 
\[
\begin{aligned}
&
a'(s)=qs^{q-1}e^{s^q},
\\
&
a''(s)=\Bigl(q(q-1)s^{q-2}+q^2s^{2q-2}\Bigr)e^{s^q},
\\
&
a'''(s)=\Bigl(q(q-1)(q-2)s^{q-3}+q^2(3q-3)s^{2q-3}+q^3s^{3q-3}\Bigr)e^{s^q},
\\
&
a''''(s)=O\Bigl(s^{4q-4}\Bigr)e^{s^q},
\end{aligned}
\]
 it follows from Lemma~\ref{lemma:3.1} that
\[
\begin{aligned}
\log F(s)
&
=
-a(s)-\log a'(s)+O\left(\frac{a''(s)}{a'(s)^2}\right)
=
-e^{s^q}-s^q-(q-1)\log s -\log q +O(e^{-s^q}),
\\
f'(s)F(s)
&
=1-\frac{a''(s)}{a'(s)^2}+\frac{3a''(s)^2-a'(s)a'''(s)}{a'(s)^4}\Bigl(1+o(1)\Bigr)
\\
&
=
1
-
\left[1+\frac{q-1}{q}\frac{1}{s^q}\right]e^{-s^q}
+
\left[
2+3\frac{q-1}{q}\frac{1}{s^q}+\frac{2q-1}{q^2}\frac{1}{s^{2q}}
\right]e^{-2s^q}\Bigl(1+o(1)\Bigr),
\\
\frac{f(s)f''(s)F(s)}{f'(s)}
&
=
1+\frac{2a''(s)^2-a'(s)a'''(s)}{a'(s)^4}
+
O\left(
\left|\frac{a''a'''}{a'^5}\right|
+
\left|\frac{a''''}{a'^4}\right|
+
\left|\frac{a''^3}{a'^6}\right|
\right)
\\
&
=
1+
\left[
1+\frac{q-1}{q}\frac{1}{s^q}+\frac{q-1}{q}\frac{1}{s^{2q}}
\right]e^{-2s^q}+O\Bigl(e^{-3s^q}\Bigr).
\end{aligned}
\]
Therefore,
\[
\frac{1}{B_1[f]}=(-\log F(s))\Bigl(1-f'(s)F(s)\Bigr)
=
1+\frac{q-1}{q}\frac{1}{s^q}+s^qe^{-s^q}+O\Bigl((\log s)e^{-s^q}\Bigr)
\to 1
\]
and
\[
\begin{aligned}
\frac{1}{B_2[f]}
&
=
f'F(-\log F)^2\left[\frac{ff''F}{f'}-1\right]
\\
&
=
1
+\frac{q-1}{q}\frac{1}{s^q}
+\frac{q-1}{q}\frac{1}{s^{2q}}
+2s^qe^{-s^q}+O\left((\log s)e^{-s^q}\right)
\\
&
\to 1
\end{aligned}
\]
as $s\to \infty$.
Furthermore,
since $\tilde u =F^{-1}[G(v)]=F^{-1}\left[\frac{1}{4}|x|^2\left(\log \frac{1}{|x|^2}+1\right)\right]$,
we obtain that
\begin{equation}\label{eq:4.44}
\begin{aligned}
&
\frac{1}{B_1[f](\tilde u)}
=
1+\frac{q-1}{q}\frac{1}{\log(-\log w)}+O\left(\frac{\log(-\log w)}{-\log w}\right),
\\
&
\frac{1}{B_2[f](\tilde u)}
=
1+\frac{q-1}{q}\frac{1}{\log(-\log w)}
+
\frac{q-1}{q}\frac{1}{\Bigl(\log(-\log w)\Bigr)^2}+O\left(\frac{\log(-\log w)}{-\log w}\right)
\end{aligned}
\end{equation}
for 
\[
w=\frac{1}{4}|x|^2\left(\log \frac{1}{|x|^2}+1\right).
\]
We now consider $B_2[g]$ for $g(s):=4\frac{e^{e^s}}{e^{2s}}=e^{e^s-2s+\log 4}$. 
Let $b(s):=e^s-2s+\log 4$. Then $b(s)$ satisfies all assumptions in Lemma~\ref{lemma:3.1}. Since
\[
b'(s)=e^s-2,\qquad b''(s)=b'''(s)=b''''(s)=e^s,
\]
we have from Lemma~\ref{lemma:3.1} that
\[
\begin{aligned}
\log G(s)
&
=-e^s+s-\log 4+O(e^{-s}),
\\
g'(s)G(s)
&
=1-e^{-s}+O(e^{-2s}),
\\
\frac{g(s)g''(s)G(s)}{g'(s)}
&
=1+\frac{2e^{2s}-(e^s-2)e^s}{(e^s-2)^4}+O(e^{-3s})
\\
&
=1+e^{-2s}+O(e^{-3s}).
\end{aligned}
\]
Therefore,
\[
\frac{1}{B_1[g]}=\frac{1}{B_2[g]}=1+O\left(se^{-s}\right)
\]
and
\[
\frac{1}{B_1[g](v)}=\frac{1}{B_2[g](v)}=1+O\left(\frac{\log(-\log w)}{-\log w}\right)
\]
for $w=\frac{1}{4}|x|^2\left(\log \frac{1}{|x|^2}+1\right)$.
All together, we obtain
\begin{equation}\label{eq:4.45}
R_1=R_2=
\left\{
\begin{aligned}
& 
O\left(\frac{1}{\log(-\log w)}\right), && q>1,
\\
&
O\left(\frac{\log(-\log w)}{-\log w}\right), && q=1
\end{aligned}
\right.
\end{equation}
for $w=\frac{1}{4}|x|^2\left(\log \frac{1}{|x|^2}+1\right)$.
This shows us that $f$ satisfies \eqref{eq:1.1a} if $q=1$, while $f$ does not satisfies \eqref{eq:1.1a} if $q>1$.

Finally let us assume that $q=1$. 
Then,
it follows from
\[
f(s)F(s)=O\left(e^{-s}\right)
\]
that
\begin{equation}\label{eq:4.46}
f\bigl(\tilde u(x)\bigr)F\bigl(\tilde u(x)\bigr)\sup_{|y|\le |x|}
\bigl(R_1(y)+R_2(y)\bigr)
=
O\left(\frac{\log(-\log w)}{(-\log  w)^{2}}\right)
\end{equation}
for $w=\frac{1}{4}|x|^2\left(\log \frac{1}{|x|^2}+1\right)$.
Therefore, \eqref{eq:example4.40}
and \eqref{eq:4.46}
yield 
the conclusion of Example~\ref{example:3.4}.
\end{proof}




\bigskip

{\noindent
{\bf Acknowledgments.}

The authors express their sincere thanks to Prof. Daisuke Naimen for informing them on the  related works \cite{DGP},\cite{GGP}.
The main part of this work was done
during a visit of the second named author
at the Department of Mathematics 
of 
Universit{\`a} degli Studi di Milano Statale.
He wishes to thank the members of the Department for their kind hospitality.
This work was partially funded by JSPS KAKENHI (grant number 19KK0349 and 21K18582) and it has been written within the activities of GNAMPA group of INDAM. 
}

%
%

\end{document}